\documentclass[lettersize,twocolumn,journal]{IEEEtran}  
\usepackage{graphicx} 
\usepackage{amsmath,bm,times} 
\usepackage{amssymb} 
\usepackage{tikz}
\usetikzlibrary{shapes,arrows,backgrounds,fit,positioning}
\usepackage{subfigure}
\usepackage{cite}
\usepackage{soul}
\usepackage{balance}
\usepackage{tabularx}
\usepackage{url}
\allowdisplaybreaks
\usepackage{multirow}
\usepackage{epstopdf}

\newtheorem{propos}{Proposition}
\newtheorem{remark}{Remark}

\newtheorem{lemma}{Lemma}

\usepackage{algorithm}
\usepackage{algpseudocode}

\newcommand{\qed}{$\hfill\blacksquare$}

\title{Semi-Explicit Solution of Some Discrete-Time Higher-Order-Cost Mean-Field-Type Control}

\author{Julian Barreiro-Gomez, \IEEEmembership{Senior Member, IEEE}, Tyrone E. Duncan  \IEEEmembership{Life Fellow Member, IEEE},\\ Bozenna Pasik-Duncan \IEEEmembership{Life Fellow Member, IEEE}, and Hamidou Tembine, \IEEEmembership{Senior Member, IEEE} 
\thanks{J. Barreiro-Gomez is with KU Center for Autonomous Robotic Systems, Department of Computer and Information Engineering, Khalifa University, Abu Dhabi 127788, UAE.(e-mail: {\tt\scriptsize julian.barreirogomez@ku.ac.ae}).}
\thanks{H. Tembine is  with Department of Electrical Engineering and Computer Science, School of Engineering  at Universit\'{e} du Qu\'{e}bec \`{a} Trois-Rivi\`{e}res UQTR, Quebec, Canada, and with Timadie, Guinaga, Grabal,  AI Mali,   TF,  WETE, MFTG, LnG Lab, CI4SI;  (e-mail: {\tt\scriptsize tembine@ieee.org}).}
\thanks{T. E. Duncan and B. Pasik-Duncan are with the Department of Mathematics, University of Kansas. Nichols Hall,
2335 Irving Hill Rd, Lawrence, KS 66045.}
\thanks{Authors gratefully acknowledge support from TIMADIE, Guinaga, Grabal, LnG Lab for the MFTG for Machine Intelligence project.
Authors gratefully acknowledge support from U.S.
Air Force Office of Scientific Research under grants
number FA9550-17-1-0259.}
\thanks{J. Barreiro-Gomez is profoundly grateful to God and to Our Lady of Lourdes for the blessings of health and life, without which his contributions to this work would not have been possible}
}

\IEEEoverridecommandlockouts

\usepackage{algorithm}
\usepackage{algpseudocode}
\usepackage[mathscr]{euscript}

\begin{document}

\maketitle

\begin{abstract}
Traditional solvable optimal control theory predominantly focuses on quadratic costs due to their analytical tractability, yet they often fail to capture critical non-linearities inherent in real-world systems including water, energy, agriculture, and financial networks. Here, we present a unified framework for solving discrete-time optimal control with higher-order state and control costs of power-law form. By building convex-completion techniques, we derive semi-explicit expressions for control laws, cost-to-go functions, and recursive coefficient dynamics across deterministic and stochastic system settings. Key contributions include variance-aware solutions under additive and multiplicative noise, extensions to mean-field-type-dependent dynamics, and conditions that ensure the positivity of recursive coefficients. In particular, we establish that higher-order costs induce less aggressive control policies compared to quadratic formulations, a finding that is validated through numerical analyses. 
\end{abstract}


\section{Introduction}
Optimal control theory has long relied on quadratic cost functions due to their analytical convenience and well-established solution techniques \cite{kirk2004optimal,stengel1994optimal}. However, many real-world control problems involve non-linear penalties that go beyond quadratic formulations, particularly in applications where deviations from a desired state incur in disproportionately large costs. Higher-order costs, such as quartic or more general power-law functions, provide a more accurate representation of such non-linearities, capturing effects like state saturation, energy constraints, and risk-sensitive behavior.
A defining characteristic of Mean-Field-Type Control (MFTC) theory is the incorporation of higher-order performance criteria, including variance, quantiles, inverse quantiles, skewness, kurtosis, and other risk measures \cite{bensoussan,Achdou_2016,Lauriere_2015}. Critically, the integration of these measures is not necessarily linear with respect to the probability distribution of individual/common state or action variables. 
Despite their practical significance, solving discrete-time optimal control problems with higher-order costs remains a challenge. The absence of closed-form solutions and the increased complexity of the resulting recursive cost structures have limited the applicability of traditional control techniques. 

In the recent literature, one may find many reported works in which semi-explicit solutions are computed to control problems both in continuous and discrete time. Also, some extensions to multi-agent systems using game theory have been introduced in the framework of mean-field-type games. Regarding semi-explicit solutions computation for risk-aware control, there are some works considering quadratic costs, state-variance, and control-variance. In  \cite{ref04}, a class of discrete-time distributed mean-variance problems is introduced and solved in a semi-explicit way. Semi-explicit solutions to continuous time linear-quadratic mean-field-type games are examined in  \cite{ref01,ref02,ref03}.  In \cite{ref0}, the authors propose a direct method to solve the game, team, and bargaining problems. As an advantage, the direct method solution approach does not require solving the Bellman-Kolmogorov equations or backward-forward stochastic differential equations of Pontryagin's type. 
The book on Mean-Field-Type Games for Engineers \cite{ref1}, presents both continuous-and discrete-time formulations of mean-field-type games, including linear-quadratic mean-field-type problems. The book discusses convex optimization techniques and equilibrium computations, and provides a MatLab toolbox for these problems \cite{Toolbox}. The work in \cite{ref2} investigates discrete-time linear-quadratic mean-field-type repeated games under different information structures, including perfect, incomplete, and imperfect information settings. The work in \cite{ref4}  employs adaptive dynamic programming techniques to solve LQ-MFTCs for unknown mean-field stochastic discrete-time systems. The study focuses on real-time learning approaches and optimal control strategies.

\begin{table*}[t!]
    \caption{Summary of Contributions Across Different Problem Types}
    \label{tab:contributions}
    \centering
    \renewcommand{\arraystretch}{1.2}
    \begin{tabular}{lccc}
        \hline
        \textbf{Problem Type} & \textbf{Optimal State} & \textbf{Optimal Control} & \textbf{Optimal Cost-to-Go}  \\
        \hline
        Deterministic Control with Higher-Order Costs & \(\checkmark\) & \(\checkmark\) & \(\checkmark\)  \\
        Variance-Aware Control with Additive Noise & \(\checkmark\) & \(\checkmark\) & \(\checkmark\) \\
        Risk-Aware Stochastic Control with Multiplicative Noise & \(\checkmark\) & \(\checkmark\) & \(\checkmark\) \\
        \hline 
    \end{tabular}
\end{table*}

On the other hand, there are some works reporting solutions to risk-aware problems considering non-quadratic costs, and higher-order moments for states and control inputs. Mean-field-type games with higher-order cost structures extend classical game theory to multi-agent systems where decision-makers optimize objectives involving non-linear or polynomial cost functions. Finding explicit solutions in these settings is challenging due to the complexity of the resulting coupled systems of partial-integro differential equations (PIDEs) or stochastic integro-differential equations (SIDEs). However, recent advances have introduced semi-explicit and fully explicit methods to solve these problems efficiently. The work in  \cite{ref5b} and \cite{ref5c} provide various solvable examples beyond the classical linear-quadratic game problems. These include quadratic-quadratic games and games with power, logarithmic, sine square, hyperbolic sine square payoffs. Non-linear state dynamics such as log-state, control-dependent regime switching, quadratic state, cotangent state, and hyperbolic cotangent state are considered. 
Hierarchical mean-field-type games with polynomial cost functions are introduced in    \cite{ref5}. In \cite{ref7} the problem of designing a collection of terminal payoff of mean-field type by interacting decision-makers to a specified terminal measure is considered. 
The work in \cite{ref6} examines a class of mean-field-type games with finite number of decision-makers with state dynamics driven by Rosenblatt processes and higher-order cost. Rosenblatt processes are non-Gaussian, non-Poisson, and non-Markov with long-range dependence.

Our results in this paper contribute to the broader understanding of non-linear control, providing a structured approach to handling cost functions of arbitrary even order. This work has implications for fields ranging from water, agriculture, energy systems, and robotics to block-chain-based finance and token economics, where risk-aware and non-linear control strategies are increasingly relevant.
In this work, we develop a semi-explicit solution framework for discrete-time optimal control settings involving non-quadratic costs. Using convex completion methods, we derive recursive formulations for optimal control strategies and cost-to-go functions, establishing conditions for their well-posedness. Our analysis extends to variance-aware control under stochastic perturbations. 
 %
%
We highlight now the contributions presented in this paper. The key contributions of this paper are summarized next and shown in Table \ref{tab:contributions}:
    
\textit{Deterministic Control with Higher-Order Costs:}
\begin{itemize}
    \item We start with a basic  framework for solving discrete-time optimal control problems with non-quadratic cost functions, specifically  the quartic terms $\bar x^{4}$ in state and quartic in control $\bar  u^{4}$.
    \item We then develop a convex-completion approach to derive semi-explicit solutions for the optimal control and the associated cost-to-go function.
    \item We establish recursive equations for cost coefficients and prove conditions ensuring their positivity, guaranteeing well-posed solutions.
    \item We extend the method to generalized even-order cost functions of the form $\bar  x^{2p}$ and $\bar  u^{2p}$.
    \item We present numerical examples illustrating the impact of increasing $p$ on optimal state trajectories, control actions and optimal cost to go.
\end{itemize}

\textit{Variance-Aware Control with Additive Noise:}
\begin{itemize}
    \item We extend the deterministic control framework to stochastic systems with additive noise, incorporating variance-aware cost functions.
    \item We then derive semi-explicit optimal control solutions for systems with quartic $var(x)+\bar x^{4}$ and $var(u)+\bar  u^{4}$, and higher-order penalties  $var(x)+\bar x^{2p}$ and $var(u)+\bar  u^{2p}$,  under stochastic perturbations.  
    \item We establish recursive conditions for cost coefficients and ensure their positivity and stability and demonstrate through semi-explicit solutions that higher-order variance-aware control leads to less aggressive control actions.
\end{itemize}

\textit{Risk-Aware Stochastic Control with Multiplicative Noise}
\begin{itemize}
    \item We generalize the variance-aware control framework to stochastic systems with multiplicative noise, capturing uncertainty-dependent dynamics.
    \item We develop recursive cost-to-go formulations accounting for both state variance and higher-order cost penalties.
    \item We derive a semi-explicit optimal control strategy that balances risk sensitivity and nonlinear cost minimization.
    \item We extend the analysis to general even-order cost functions and solve it semi-explicitly.
\end{itemize}
%
%
\vspace{-0.3cm}
\subsection{How does the Convex Completion Works?}
As it was presented before in the literature review, there are many works on the computation of semi-explicit solutions to mean-field-type control problems by following the so-called direct method. One of the key elements in this method, in the framework of linear-quadratic settings, is the square completion in order to optimize the Hamiltonian over the control actions. In this work, we do not longer work with the square completion as the proposed costs in this paper go beyond the quadratic structure, and examines higher order costs. Instead, this work uses convex completion as briefly introduced next.

Let $X$ be a real vector space and $X^*$ be its dual space. Let $f: X \to \mathbb{R} \cup \{+\infty\}$ be a proper, convex, and lower semi-continuous function. The Legendre-Fenchel conjugate\footnote{Also known as convex conjugate.} (or transform) of $f$ is defined as:
$$f^*(y) = \sup_{x \in X} \{\langle y, x \rangle - f(x)\},$$
where $\langle y, x \rangle$ denotes the dual pairing between $y \in X^*$ and $x \in X$.
The convexity inequality, also known as Fenchel's inequality or Young's inequality, states that for any $x \in X$ and $y \in X^*$:
$$f(x) + f^*(y) \geq \langle y, x \rangle.$$
The convex completion method uses this inequality as follows. After making the difference between the candidate cost and control cost functional, one gets a function of the control, say $f(u).$ We complete this cost by adding  and removing $f^*(v) - \langle v, u \rangle.$
$$ f(u)+f^*(v) - \langle v, u \rangle \geq 0 $$ for control actions $(u,v).$
The equality holds if and only if $v \in \partial f(u)$, where $\partial f(u)$ is the sub-differential of $f$ at $u$. The sub-differential is a generalization of the gradient for non-differentiable functions. 
This extra term added to the difference cost goes to the coefficient matching terms making the solution semi-explicit.

\vspace{-0.3cm}
\subsection{Convexity of higher-order terms}

This paper is going to consider higher order terms, i.e., beyond the quadratic case, with a particular structure. Next, Lemma \ref{lema:highorderterms} shows a convexity property that the development of the paper relies on.
\begin{lemma}
\label{lema:highorderterms}
Let $p\geq 1, a\neq 0, b\neq 0.$ Then the mapping $z \mapsto f(z) = z^{2p} + (az+b)^{2p}$ is strictly convex. \hfill $\square$
\end{lemma}

\vspace{-0.3cm}
\subsection{Structure of the paper}
The rest of the article is organized as follows. 
Deterministic control problems with arbitrary order are presented in  Section 
 \ref{sec:label:1}. Section \ref{sec:label:2} discusses variance-aware control with additive noise and higher-order cost.  
   Section \ref{sec:label:3} presents  risk-aware stochastic control with multiplicative noise. Section \ref{sec:Multiplicative_mean_field_noise} presents the case with multiplicative and both state-control and mean-field-type noise. Section \ref{sec:proofs} presents the proofs corresponding to the results shown in this paper. Numerical examples are introduced in Section \ref{sec:examples}. Finally, concluding remarks are drawn in Section \ref{sec:conclusions}.

%
%
\section{Deterministic control with higher-order costs} \label{sec:label:1}

In this section, we explicitly solve the discrete-time optimal control problem with higher-order costs. We derive closed-form solutions for the optimal control, the optimal cost, and the recursive coefficients for the cost-to-go function. Additionally, we extend the methodology to address general higher-order costs of the form $\bar{u}^{2p}$ for arbitrary $p \geq 1$, where $\bar{u}$ denotes the expected control input.
\subsection{Optimal Control with Quartic Cost}
We consider the following discrete-time dynamical system:
\begin{equation} \label{eqstate000}
    \bar{x}_{k+1} = \bar{a}_k \bar{x}_k + \bar{b}_k \bar{u}_k,
\end{equation}
with $k \in \{0, \dots, N-1\}$, where $\bar{x}_k \in \mathbb{R}$ is the system state, $\bar{u}_k \in \mathbb{R}$ is the control input, and $\bar{a}_k, \bar{b}_k \in \mathbb{R} $ are system parameters. The dynamics in \eqref{eqstate000} evolves within a discrete-time horizon $[0..N]:= [0,N] \cap \mathbb{Z}$ with $N \in \mathbb{Z}_+$. The initial state $\bar{x}_0$ is given. The goal is to minimize the cost function:
\begin{equation}  \label{eqcost000}
    L(\bar{x}_0,\bar{u}) = \bar{q}_N \bar{x}_N^4 + \sum_{k=0}^{N-1} \left(\bar{q}_k \bar{x}_k^4 + \bar{r}_k \bar{u}_k^4\right),
\end{equation}
where $\bar{q}_k, \bar{q}_N, \bar{r}_k > 0$ are weights penalizing the state and control deviations. Therefore, the optimal control problem is
\begin{align}
\label{eq:problem1}
    \min_{\bar{u} \in \mathcal{\bar U}_k} L(\bar{x}_0,\bar{u}),~\mathrm{s.t.}~\eqref{eqstate000},~\bar{x}_0~\mathrm{given}.
\end{align}
In order to solve the problem in \eqref{eq:problem1} 
we need to specify the information structure and control action. The set of control actions  at time step $k$ is 
$\bar{U}_k= \{ \bar{u}_k\in \mathbb{R} | \bar{u}_k^4 < +\infty \}=\mathbb{R},$  
and the set 
\begin{align}
    \label{eq:cal_bar_U}
    \mathcal{\bar U}_k= \left\{  (\bar{u}_k)_{k\in \{0, \dots, N-1\}} \ \bigg| \begin{array}{l}
   \bar{u}_k(.) \in \bar{U}_k , \\
    \bar{u}_k~  \mbox{Lebesgue measurable} 
\end{array}  \hspace{-0.1cm}\right\}.
\end{align}
We restrict our attention to state-feedback control, i.e.,  
$\bar{u}_k( \bar{x}_k).$ This means that instead of considering the entire past history  $(  \bar{x}_0,  \bar{x}_1, \ldots,  \bar{x}_{k-1},  \bar{x}_{k})$, we choose to work only with the current state $ \bar{x}_k.$ 

Note that for any $\bar u\in \mathcal{\bar U},$ the system  (\ref{eqstate000})  admits unique solution denoted by $ \bar{x}:=  \bar{x}^{ \bar{x}_0,\bar u}$ meaning it is determined by the initial state and the control strategy employed over the path. In this particular case the solution is explicitly obtained by  recursion  of linear (affine) operators as shown next in Proposition \ref{propos:problem1}. 

\vspace{0.1cm}
\begin{propos}
\label{propos:problem1}
    The optimal control input for the Problem in \eqref{eq:problem1} 
    with system dynamics in \eqref{eqstate000} and cost functional \eqref{eqcost000}
    is
    \begin{equation}
    \label{eq:u_optimal1}
    \bar{u}_k^* = -\frac{\left(\frac{\bar{\alpha}_{k+1} \bar{b}_k}{\bar{r}_k}\right)^{1/3}\bar{a}_k }{1 + \left(\frac{\bar{\alpha}_{k+1} \bar{b}_k}{\bar{r}_k}\right)^{1/3} \bar{b}_k} \bar{x}_k,
    \end{equation}
where $\bar{\alpha}$ solves the following backward recursive equation:
\begin{align}
\label{eq:alpha_optimal1}
    \bar{\alpha}_k = \bar{q}_k + \bar{r}_k &\left(\frac{\left(\frac{\bar{\alpha}_{k+1} \bar{b}_k}{\bar{r}_k}\right)^{1/3} \bar{a}_k }{1 + \left(\frac{\bar{\alpha}_{k+1} \bar{b}_k}{\bar{r}_k}\right)^{1/3} \bar{b}_k}\right)^4 \notag\\
    &+ \bar{\alpha}_{k+1} \left(\bar{a}_k -  \frac{\left(\frac{\bar{\alpha}_{k+1} \bar{b}_k}{\bar{r}_k}\right)^{1/3} \bar{a}_k \bar{b}_k}{1 + \left(\frac{\bar{\alpha}_{k+1} \bar{b}_k}{\bar{r}_k}\right)^{1/3} \bar{b}_k}\right)^4,
\end{align}
with terminal condition $\bar{\alpha}_N= \bar{q}_N$. \hfill $\square$
\end{propos}

Considering quartic cost functional implies that the cost can only get non-negative values as $\bar{q}_k, \bar{q}_N, \bar{r}_k > 0$. Therefore, according to the ansatz for the value function that is used in the proof (see Section \ref{sec:proofs}), the recursive solution for $\bar{\alpha}_k$ in \eqref{eq:alpha_optimal1} should be positive. The Remark \ref{rmk:coeffs} below shows the conditions under which there is existence of positive $\bar{\alpha}_k$.

\vspace{0.1cm}
\begin{remark}
\label{rmk:coeffs}
    The coefficients $\bar{\alpha}_k$ remain positive in \eqref{eq:alpha_optimal1} if:
\begin{itemize}
    \item The cost parameters $\bar{q}_k, \bar{q}_N, \bar{r}_k > 0,$
    \item The system parameters $\bar{a}_k, \bar{b}_k$ are bounded,
    \item The denominator in $\bar{u}_k^*$ satisfies  $1 + \left(\frac{\bar{\alpha}_{k+1} \bar{b}_k}{\bar{r}_k}\right)^{1/3} \bar{b}_k > 0.$
\end{itemize}
These conditions are trivially satisfied in most systems, ensuring $\bar{\alpha}_k > 0$, for all  $k \in \{0, 1, \ldots, N \}$. \hfill $\square$
\end{remark}

\vspace{0.1cm}
In the next Section, we present how to generalize the power of the system state in the cost functional to be any even power.

\subsection{Extension to $ \bar{x}^{2p}$ Costs}

We now address discrete-time optimal control problems with costs of the form $\bar{x}_k^{2p}$ and $\bar{u}_k^{2p}$, where $p \geq 1. $ The solutions include the optimal control law, the recursive relation for the cost-to-go coefficient $\bar{\alpha}_k,$ and the expression for the optimal cost. We consider the same discrete-time system as in \eqref{eqstate000}. 
The objective is to minimize the $2p$-order cost function:
\begin{align}
\label{eq:cost_2p}
    L(x_0,\bar{u}) = \bar{q}_N \bar{x}_N^{2p} + \sum_{k=0}^{N-1} \left(\bar{q}_k \bar{x}_k^{2p} + \bar{r}_k \bar{u}_k^{2p}\right),
\end{align}
where $\bar{q}_k, \bar{q}_N, \bar{r}_k > 0$  are weights in the cost functional. \textcolor{black}{The new set of control actions at time step $k$ is 
$\bar{U}_k= \{ \bar{u}_k\in \mathbb{R} | \bar{u}_k^{2p} < +\infty \}=\mathbb{R},$ with $\mathcal{\bar U}_k$ as in \eqref{eq:cal_bar_U}.} The new result is presented in Proposition \ref{propos:problem2} below.
\begin{propos}
\label{propos:problem2}
    The optimal control input for the Problem in \eqref{eq:problem1} with system dynamics in \eqref{eqstate000} and general power cost functional \eqref{eq:cost_2p} is given by 
    \begin{align}
    \label{eq:u_optimal2}
    \bar{u}^*_k = -\frac{\left(\frac{\bar{\alpha}_{k+1} \bar{b}_k}{\bar{r}_k}\right)^{1/(2p-1)} \bar{a}_k }{1 + \left(\frac{\bar{\alpha}_{k+1} \bar{b}_k}{\bar{r}_k}\right)^{1/(2p-1)} \bar{b}_k}\bar{x}_k,
\end{align}
where $\bar{\alpha}$ solves the following backward recursive equation:
\begin{align}
\label{eq:alpha_optimal2}
    \bar{\alpha}_k &= \bar{q}_k + \bar{r}_k \left(\frac{\left(\frac{\bar{\alpha}_{k+1} \bar{b}_k}{\bar{r}_k}\right)^{1/(2p-1)} \bar{a}_k}{1 + \left(\frac{\bar{\alpha}_{k+1} \bar{b}_k}{\bar{r}_k}\right)^{1/(2p-1)} \bar{b}_k}\right)^{2p} \notag\\
    &\quad + \bar{\alpha}_{k+1} \left(\bar{a}_k - \bar{b}_k \frac{\left(\frac{\bar{\alpha}_{k+1} \bar{b}_k}{\bar{r}_k}\right)^{1/(2p-1)} \bar{a}_k}{1 + \left(\frac{\bar{\alpha}_{k+1} \bar{b}_k}{\bar{r}_k}\right)^{1/(2p-1)} \bar{b}_k}\right)^{2p},
\end{align}
with terminal condition $\bar{\alpha}_N= \bar{q}_N$ depending on the cost weight. \hfill $\square$
\end{propos}

The result presented in Proposition \ref{propos:problem2} is more general than the quadratic cost that has been reported in the literature. Hence, one can easily retrieve the quadratic solution by adjusting the $p$-value as shown in the Remark \ref{rmk:remark2} below.

\vspace{0.1cm}
\begin{remark}
\label{rmk:remark2}
For $p=1$ the latter equation becomes
\begin{align*}
    \bar{\alpha}_k &= \bar{q}_k + \frac{\bar{\alpha}_{k+1} \bar{a}_k^2 \bar{r}_k}{\bar{r}_k+ \bar{\alpha}_{k+1} \bar{b}_k^2},&
    \bar{\alpha}_N &= \bar{q}_N,
\end{align*}
which is the Riccati equation for the discrete-time Linear-Quadratic Regulator (LQR) case \cite{kirk2004optimal,stengel1994optimal}.  \hfill $\square$
\end{remark}

\section{Variance-aware control with additive noise}  \label{sec:label:2}
Non-quadratic costs, such as $ \bar{x}^4$ or $\bar{u}^4$, with $\bar{x} = \mathbb{E}[x]$ and $\bar{u} = \mathbb{E}[u]$, introduce another class of  realistic representation of these non-linearities, emphasizing large deviations in state and control.
We  solve the variance-aware discrete-time stochastic optimal control problem with higher-order costs. We derive semi-explicit solutions for the  variance-aware stochastic optimal control, the optimal  variance-aware cost, and the recursive coefficients for the  variance-aware cost-to-go function. Additionally, we extend the methodology to address general higher-order costs $\bar{x}^{2p}$ for arbitrary $p \geq 1.$
\subsection{Variance-aware Quartic Cost with additive noise}
We consider the following discrete-time dynamical system:
\begin{align}
\label{eq:dynamics2}
    {x}_{k+1} = \bar{a}_k {x}_k +  \bar{b}_k {u}_k +\epsilon_{k+1} ,
\end{align}
with $k \in \{ 0, \dots, N-1\}$, where $x_k \in \mathbb{R}$ is the system state, $u_k \in \mathbb{R}$ is the control input, and $ \bar{a}_k,  \bar{b}_k \in \mathbb{R} $ are system 
parameters. The initial state $x_0$ is a given random variable that is independent of $\epsilon$.  Also, $ \{\epsilon_{k}\}_k$ is a zero-mean random process, with $\epsilon_0=0$.  
The state dynamics in \eqref{eq:dynamics2} can be re-written as 
\begin{align*}
   {x}_{k+1} =   \bar{a}_k \bar{x}_k + \bar{b}_k \bar{u}_k+
   \bar{a}_k (x_k-\bar{x}_k) + \bar{b}_k (u_k-\bar{u}_k)+\epsilon_{k+1},
\end{align*}
with $k \in \{ 0, \dots, N-1\}$. The goal is to minimize the expected value of the following cost function:
\begin{align*}
    L(x_0,u) &=   {q}_k {var(x_N)} +\bar{q}_N  \bar{x}_N^4 \\
    &+ \sum_{k=0}^{N-1} \left( {q}_k {var(x_k)}+\bar{q}_k  \bar{x}_k^4 + {r}_k {var(u_k)}+ \bar{r}_k  \bar{u}_k^4\right),
\end{align*}
where $  q_k,{q}_N, r_k, \bar{q}_k,  \bar{q}_N,   \bar{r}_k > 0$ are weights penalizing the system state and control deviations. The variance of the system state and control input are given by $var(x) = \mathbb{E}[(x-\bar{x})^2]$, and $var(u) = \mathbb{E}[(u-\bar{u})^2]$, respectively. Therefore, the variance-aware control problem becomes
\begin{align}
\label{eq:problem2}
    \min_{u \in \mathcal{U}_k} \mathbb{E} [L(x_0,u)],~\mathrm{s.t.}~\eqref{eq:dynamics2},~x_0~\mathrm{given}.
\end{align}
\textcolor{black}{The set of control actions  at time step $k$ is 
${U}_k= \{ {u}_k\in \mathbb{R} | \bar{u}_k^4 < +\infty \}=\mathbb{R},$  
and the set 
\begin{align}
\label{eq:cal_U}
\mathcal{U}_k=  \left\{  ({u}_k)_{k\in \{0, \dots, N-1\}} \ \hspace{-0.1cm} \bigg| \hspace{-0.1cm} \begin{array}{l}
   {u}_k(.) \in {U}_k , \\
    {u}_k~  \mbox{Lebesgue measurable} 
\end{array} \hspace{-0.1cm} \right\}.    
\end{align}} 
The solution is presented in Proposition \ref{propos:proposition_3} below. 

\begin{propos}
\label{propos:proposition_3}
    The optimal control input for the variance-aware problem in \eqref{eq:problem2} is given by
    \begin{equation}
    \label{eq:optimal_u_3}
    {u}_k^* = 
-  \frac{ {\alpha}_{k+1} \bar{b}_k   \bar{a}_k }
{r_k +\alpha_{k+1} \bar{b}^2_k} (x_k-\bar{x}_k)
 -\frac{\left(\frac{\bar{\alpha}_{k+1} \bar{b}_k}{\bar{r}_k}\right)^{1/3} \bar{a}_k }{1 + \left(\frac{\bar{\alpha}_{k+1} \bar{b}_k}{\bar{r}_k}\right)^{1/3} \bar{b}_k} \bar{x}_k,
\end{equation}
where $\bar{\alpha}$, $\bar{\gamma}$, and ${\alpha}$ solve the following recursive equations:
\begin{subequations}
\label{eq:recursive_2}
\begin{align}
    \bar{\alpha}_k &= \bar{q}_k + \bar{r}_k \left(\frac{\left(\frac{\bar{\alpha}_{k+1} \bar{b}_k}{\bar{r}_k}\right)^{1/3} \bar{a}_k }{1 + \left(\frac{\bar{\alpha}_{k+1} \bar{b}_k}{\bar{r}_k}\right)^{1/3} \bar{b}_k}\right)^4 \notag\\
    &\quad\quad\quad\quad\, + \bar{\alpha}_{k+1} \left(\bar{a}_{k}-  \frac{\left(\frac{\bar{\alpha}_{k+1} \bar{b}_k}{\bar{r}_k}\right)^{1/3} \bar{a}_k \bar{b}_k}{1 + \left(\frac{\bar{\alpha}_{k+1} \bar{b}_k}{\bar{r}_k}\right)^{1/3} \bar{b}_k}\right)^4, \\
       \bar{\gamma}_{k} &= \bar{\gamma}_{k+1} +  {\alpha}_{k+1} \mathbb{E}\left(  \epsilon_{k+1}^2  \right), \\
         {\alpha}_k &= {q}_k + {r}_k 
         \left(
          \frac{ {\alpha}_{k+1} \bar{b}_k   \bar{a}_k }
{r_k +\alpha_{k+1} \bar{b}^2_k}
         \right)^2 \notag\\
         &\quad\quad\quad\quad\quad\quad\quad\, + {\alpha}_{k+1}
          \left(\bar{a}_{k}-  
           \frac{ {\alpha}_{k+1} \bar{b}^2_k   \bar{a}_k }
{r_k +\alpha_{k+1} \bar{b}^2_k}  
          \right)^2. 
\end{align}
\end{subequations}
with $\bar{\alpha}_N = \bar{q}_N$, $\bar{\gamma}_N =0$, and ${\alpha}_N = {q}_N$. \hfill $\square$
\end{propos}

Given the structure of the cost functional and its positiveness, the coefficients for the ansatz should also satisfy positiveness as shown next in Remark \ref{rmk:remark3}.   
\begin{remark}
\label{rmk:remark3}
The coefficients ${\alpha}_k,  \bar{\alpha}_k$ remain positive if:
\begin{enumerate}
    \item $ {q}_k, {q}_N, {r}_k, \bar{q}_k, \bar{q}_N, \bar{r}_k > 0,$
    \item $ \bar{a}_k, \bar{b}_k$ are bounded,
    \item The denominator in $ \bar{u}_k^*$ satisfies  $1 + \left(\frac{\bar{\alpha}_{k+1} \bar{b}_k}{\bar{r}_k}\right)^{1/3} \bar{b}_k > 0.$
    \item The denominator in $ u_k^*-\bar{u}_k^*$ satisfies  $ (r_k +\alpha_{k+1} \bar{b}^2_k)>0.$
\end{enumerate}
These conditions are trivially satisfied in most systems, ensuring $\alpha_k > 0$  and  $ \bar{\alpha}_k>0$ for all  $k \in \{0, 1, \ldots, N \}$. \hfill $\square$
\end{remark}

In the coming Section, we present a generalization of the power to any even high-order cost.

\subsection{Extension to $ \bar{x}^{2p}$ Costs}

We now address discrete-time optimal control problems with costs of the form $x_k^{2p}$ and $u_k^{2p}$, where  $p \geq 1. $ The solutions include the optimal control law, the recursive relation for the cost-to-go coefficient $\alpha_k,$ and the expression for the optimal cost.
We consider the same discrete-time system as in \eqref{eq:dynamics2}. The objective is to minimize the variance-aware $2p$-order cost function:
\begin{align}
\label{eq:cost_2b}
    L(x_0,u) &=  {q}_N {var(x_N)} +
    \bar{q}_N \bar{x}_N^{2p} \\
    &+ \sum_{k=0}^{N-1} \left({q}_k {var(x_k)} +
    \bar{q}_k \bar{x}_k^{2p} + {r}_k {var(u_k)}+
     \bar{r}_k \bar{u}_k^{2p}\right), \notag
\end{align}
where ${q}_k, {q}_N, {r}_k , \bar{q}_k, \bar{q}_N, \bar{r}_k > 0$  are weights. \textcolor{black}{The new set of control actions at time step $k$ is 
${U}_k= \{ {u}_k\in \mathbb{R} | \bar{u}_k^{2p} < +\infty \}=\mathbb{R},$ with $\mathcal{U}_k$ as in \eqref{eq:cal_U}.}
\begin{propos}
\label{propos:proposition_4}
    The optimal control input for the variance-aware problem in \eqref{eq:problem2} with dynamics in \eqref{eq:dynamics2} and cost functional \eqref{eq:cost_2b}, is given by
    \begin{align}
    \label{eq:optimal_u_4}
    {u}^*_k &=  -  \frac{ {\alpha}_{k+1} \bar{b}_k   \bar{a}_k }
{r_k +\alpha_{k+1} \bar{b}^2_k} (x_k-\bar{x}_k) \notag\\
  &\quad\quad\quad\quad\quad\quad\quad\quad -\frac{\left(\frac{\bar{\alpha}_{k+1} \bar{b}_k}{\bar{r}_k}\right)^{1/(2p-1)} \bar{a}_k }{1 + \left(\frac{\bar{\alpha}_{k+1} \bar{b}_k}{\bar{r}_k}\right)^{1/(2p-1)} \bar{b}_k} \bar{x}_k.
\end{align}
where $\bar{\alpha}$, $\bar{\gamma}$, and ${\alpha}$ solve the following recursive equations:
\begin{subequations}
\begin{align}
    \bar{\alpha}_k &= \bar{q}_k + \bar{r}_k \left(\frac{\left(\frac{\bar{\alpha}_{k+1} \bar{b}_k}{\bar{r}_k}\right)^{1/(2p-1)} \bar{a}_k}{1 + \left(\frac{\bar{\alpha}_{k+1} \bar{b}_k}{\bar{r}_k}\right)^{1/(2p-1)} \bar{b}_k}\right)^{2p} \notag\\
    &+ \bar{\alpha}_{k+1} \left(\bar{a}_{k}- \bar{b}_k \frac{\left(\frac{\bar{\alpha}_{k+1} \bar{b}_k}{\bar{r}_k}\right)^{1/(2p-1)} \bar{a}_k}{1 + \left(\frac{\bar{\alpha}_{k+1} \bar{b}_k}{\bar{r}_k}\right)^{1/(2p-1)} \bar{b}_k}\right)^{2p},
     \\
      \bar{\gamma}_{k} &= \bar{\gamma}_{k+1} +  {\alpha}_{k+1} \mathbb{E}\left( \epsilon_{k+1}^2  \right), \\
        {\alpha}_k &= {q}_k + {r}_k 
         \left(
          \frac{ {\alpha}_{k+1} \bar{b}_k   \bar{a}_k }
{r_k +\alpha_{k+1} \bar{b}^2_k}
         \right)^2 \notag\\
         &\quad\quad\quad\quad\quad\quad + {\alpha}_{k+1}
          \left(\bar{a}_{k}-  
           \frac{ {\alpha}_{k+1} \bar{b}^2_k   \bar{a}_k }
{r_k +\alpha_{k+1} \bar{b}^2_k}  
          \right)^2, 
\end{align}
\end{subequations}
with $\bar{\alpha}_N = \bar{q}_N$, $\bar{\gamma}_N =0$, and ${\alpha}_N = {q}_N$. \hfill $\square$
\end{propos}

Interestingly, by comparing the risk-aware result in Proposition \ref{propos:proposition_4} with the risk-free result in Proposition \ref{propos:problem2}, we can appreciate that the obtained structure for the recursive equation related to the deterministic components are the same, and that there are some emerging recursive equations linked to the variance minimization.

\begin{remark}
Note that for $p=1$ the latter equation becomes
\begin{align*}
    \bar{\alpha}_k &= \bar{q}_k + \frac{\bar{\alpha}_{k+1} \bar{a}_k^2 \bar{r}_k}{\bar{r}_k+ \bar{\alpha}_{k+1} \bar{b}_k^2},&
      \bar{\alpha}_N&= \bar{q}_N,
\end{align*} 
which  is the Riccati equation for discrete-time Linear-Quadratic Regulator (LQR) case, and 
\begin{align*}
    f_0(x_0, \bar{x}_0) = {\alpha}_0 var({x}_0) + \bar{\alpha}_0 \bar{x}_0^{2p}+ \bar\gamma_0.
\end{align*} 
is the optimal cost. \hfill $\square$
\end{remark}
In the following Section we combine the consideration of high-order costs and the variance minimization. This is a risk-aware problem using the second moment.

\section{Variance-Aware Stochastic Cost with Multiplicative Noise}  \label{sec:label:3}

We consider the following discrete-time dynamical system with multiplicative state-and-mean state dependent noise:
\begin{align}
\label{eq:dynamics_3}
    {x}_{k+1} = \bar{a}_k {x}_k +  \bar{b}_k {u}_k + (x_k-\bar{x}_k)\epsilon_{k+1} ,
\end{align}
with $k \in \{ 0, \dots, N-1\}$, where $x_k \in \mathbb{R}$ is the state, $u_k \in \mathbb{R}$ is the control input, and $ \bar{a}_k,  \bar{b}_k \in \mathbb{R} $ are system 
parameters. The initial state $x_0$ is a given random variable that is independent of $\epsilon.$  $ \{\epsilon_{k}\}_k$ is a zero-mean random process, with $\epsilon_0=0.$ 
The state dynamics can be re-written as 
\begin{align*}
   {x}_{k+1} &=   \bar{a}_k \bar{x}_k + \bar{b}_k \bar{u}_k+
   \bar{a}_k (x_k-\bar{x}_k) \\
   &\quad\quad\quad\quad\quad\quad\quad\quad + \bar{b}_k (u_k-\bar{u}_k)+ (x_k-\bar{x}_k)\epsilon_{k+1},
\end{align*}
with $k \in \{ 0, \dots, N-1\}$. 
The goal is to minimize the expected value of the cost function:
\begin{align}
\label{eq:cost_multiplicative}
    L(x_0,u) &=   {q}_k {var(x_N)} +\bar{q}_N  \bar{x}_N^4 \\
    &+ \sum_{k=0}^{N-1} \left( {q}_k {var(x_k)} +\bar{q}_k  \bar{x}_k^4 + {r}_k {var(u_k)} + \bar{r}_k  \bar{u}_k^4\right),\notag
\end{align}
where $q_k,{q}_N, r_k, \bar{q}_k,  \bar{q}_N,   \bar{r}_k > 0$ are weights penalizing the state and control deviations. 
Therefore, the control problem is
\begin{align}
\label{eq:problem3}
    \min_{u \in \mathcal{U}_k} \mathbb{E} [L(x_0,{u})],~\mathrm{s.t.}~\eqref{eq:dynamics_3},~x_0~\mathrm{given}.
\end{align}
\textcolor{black}{The set of control actions at time step $k$ is 
${U}_k= \{ {u}_k\in \mathbb{R} | \bar{u}_k^{4} < +\infty \}=\mathbb{R},$ with $\mathcal{U}_k$ as in \eqref{eq:cal_U}.}
\begin{propos}
\label{propos:proposition_5}
    The optimal control input for the variance-aware problem with multiplicative noise in \eqref{eq:problem3} with dynamics in \eqref{eq:dynamics_3} and cost functional \eqref{eq:cost_multiplicative}, is given by
    \begin{align}
    \label{eq:optimal_u_5}
    {u}_k^* &= 
-  \frac{ {\alpha}_{k+1} \bar{b}_k   \bar{a}_k }
{r_k +\alpha_{k+1} \bar{b}^2_k} (x_k-\bar{x}_k) \notag\\
 &\quad\quad\quad\quad\quad\quad\quad\quad\quad\quad -\frac{\left(\frac{\bar{\alpha}_{k+1} \bar{b}_k}{\bar{r}_k}\right)^{1/3} \bar{a}_k }{1 + \left(\frac{\bar{\alpha}_{k+1} \bar{b}_k}{\bar{r}_k}\right)^{1/3} \bar{b}_k} \bar{x}_k.
\end{align}
where $\alpha$, and $\bar{\alpha}$ solve the following recursive equations:
\begin{subequations}
\label{eq:resursive_eqs}
\begin{align}
    \bar{\alpha}_k &= \bar{q}_k + \bar{r}_k \left(\frac{\left(\frac{\bar{\alpha}_{k+1} \bar{b}_k}{\bar{r}_k}\right)^{1/3} \bar{a}_k }{1 + \left(\frac{\bar{\alpha}_{k+1} \bar{b}_k}{\bar{r}_k}\right)^{1/3} \bar{b}_k}\right)^4 \notag\\
    &\quad\quad\quad\quad+ \bar{\alpha}_{k+1} \left(\bar{a}_{k}-  \frac{\left(\frac{\bar{\alpha}_{k+1} \bar{b}_k}{\bar{r}_k}\right)^{1/3} \bar{a}_k \bar{b}_k}{1 + \left(\frac{\bar{\alpha}_{k+1} \bar{b}_k}{\bar{r}_k}\right)^{1/3} \bar{b}_k}\right)^4, \\
         {\alpha}_k &= {q}_k + {r}_k 
         \left(
          \frac{ {\alpha}_{k+1} \bar{b}_k   \bar{a}_k }
{r_k +\alpha_{k+1} \bar{b}^2_k}
         \right)^2 \notag\\
         &+ {\alpha}_{k+1}
          \left(\bar{a}_{k}-  
           \frac{ {\alpha}_{k+1} \bar{b}^2_k   \bar{a}_k }
{r_k +\alpha_{k+1} \bar{b}^2_k}  
          \right)^2  +  {\alpha}_{k+1} \mathbb{E}\left( \epsilon_{k+1}^2  \right), 
\end{align}
\end{subequations}
with $\bar{\alpha}_N= \bar{q}_N$, and ${\alpha}_N= {q}_N$ as terminal conditions. \hfill $\square$
\end{propos}

\begin{remark}
The coefficients ${\alpha}_k,  \bar{\alpha}_k$ remain positive if:
\begin{enumerate}
    \item $ {q}_k, {q}_N, {r}_k, \bar{q}_k, \bar{q}_N, \bar{r}_k > 0,$
    \item $ \bar{a}_k, \bar{b}_k$ are bounded,
    \item The denominator in $ \bar{u}_k^*$ satisfies  $1 + \left(\frac{\bar{\alpha}_{k+1} \bar{b}_k}{\bar{r}_k}\right)^{1/3} \bar{b}_k > 0.$
    \item The denominator in $ u_k-\bar{u}_k^*$ satisfies  $ (r_k +\alpha_{k+1} \bar{b}^2_k)>0.$ 
\end{enumerate}    
These conditions are trivially satisfied in most systems, ensuring $\alpha_k > 0$  and  $ \bar{\alpha}_k$ for all  $k \in \{0, 1, \ldots, N \}$. \hfill $\square$
\end{remark}

Similar to what we have presented in the previous case, we now intend to generalize the power for the terms cost functional that is independent of the variance.

\subsection{Extension to $ \bar{x}^{2p}$ Costs}

We now provide for discrete-time optimal control problems with costs of the form $x_k^{2p}$ and $u_k^{2p}$, where  $p \geq 1. $ The solutions include the optimal control law, the recursive relation for the cost-to-go coefficient $\alpha_k,$ and the expression for the optimal cost. We consider the same discrete-time system as in \eqref{eq:dynamics_3}. 
%
Now, the objective is to minimize the variance-aware 2p-order cost function:
\begin{align}
\label{eq:var2p}
    L(x_0,u) &=  {q}_N {var(x_N)} +
    \bar{q}_N \bar{x}_N^{2p} \\
    &+ \sum_{k=0}^{N-1} \left(  {q}_k {var(x_k)} +
    \bar{q}_k \bar{x}_k^{2p} + {r}_k {var(u_k)}+
     \bar{r}_k \bar{u}_k^{2p}\right),\notag
\end{align}
where ${q}_k, {q}_N, {r}_k , \bar{q}_k, \bar{q}_N, \bar{r}_k > 0$  are weights. \textcolor{black}{The set of control actions at time step $k$ is 
${U}_k= \{ {u}_k\in \mathbb{R} | \bar{u}_k^{2p} < +\infty \}=\mathbb{R},$ with $\mathcal{U}_k$ as in \eqref{eq:cal_U}.}
\begin{propos}
\label{propos:proposition_6}
    The optimal control input for the variance-aware problem with multiplicative noise in \eqref{eq:problem3} with dynamics in \eqref{eq:dynamics_3} and cost functional \eqref{eq:var2p}, is given by
\begin{align}
    \label{eq:optimal_u_6}
    {u}^*_k &=   - \frac{ {\alpha}_{k+1} \bar{b}_k   \bar{a}_k }
{r_k +\alpha_{k+1} \bar{b}^2_k} (x_k-\bar{x}_k) \notag\\
&\quad\quad\quad\quad\quad\quad\quad\quad -\frac{\left(\frac{\bar{\alpha}_{k+1} \bar{b}_k}{\bar{r}_k}\right)^{1/(2p-1)} \bar{a}_k }{1 + \left(\frac{\bar{\alpha}_{k+1} \bar{b}_k}{\bar{r}_k}\right)^{1/(2p-1)} \bar{b}_k} \bar{x}_k,
\end{align}
where ${\alpha}$, and $\bar{\alpha}$ solve the following recursive equations:
\begin{subequations}
\label{eq:resursive_eqs_2}
\begin{align}
\bar{\alpha}_k &= \bar{q}_k + \bar{r}_k \left(\frac{\left(\frac{\bar{\alpha}_{k+1} \bar{b}_k}{\bar{r}_k}\right)^{1/(2p-1)} \bar{a}_k}{1 + \left(\frac{\bar{\alpha}_{k+1} \bar{b}_k}{\bar{r}_k}\right)^{1/(2p-1)} \bar{b}_k}\right)^{2p} \notag\\
&+ \bar{\alpha}_{k+1} \left(\bar{a}_{k}- \bar{b}_k \frac{\left(\frac{\bar{\alpha}_{k+1} \bar{b}_k}{\bar{r}_k}\right)^{1/(2p-1)} \bar{a}_k}{1 + \left(\frac{\bar{\alpha}_{k+1} \bar{b}_k}{\bar{r}_k}\right)^{1/(2p-1)} \bar{b}_k}\right)^{2p},
     \\
        {\alpha}_k &= {q}_k + {r}_k 
         \left(
          \frac{ {\alpha}_{k+1} \bar{b}_k   \bar{a}_k }
{r_k +\alpha_{k+1} \bar{b}^2_k}
         \right)^2 \notag\\
         &+ {\alpha}_{k+1}
          \left(\bar{a}_{k}-  
           \frac{ {\alpha}_{k+1} \bar{b}^2_k   \bar{a}_k }
{r_k +\alpha_{k+1} \bar{b}^2_k}  
          \right)^2 +   {\alpha}_{k+1} \mathbb{E}\left( \epsilon_{k+1}^2  \right), 
\end{align}
\end{subequations}
with $\bar{\alpha}_N= \bar{q}_N$, and ${\alpha}_N = {q}_N$ as terminal conditions. \hfill $\square$
\end{propos}

Note that, when comparing additive noise case in Proposition \ref{propos:proposition_4} with the multiplicative noise case in Proposition \ref{propos:proposition_6}, we observe an equivalent recursive equation for the deterministic part $\bar{\alpha}_k$; however, the recursive equations in $\alpha_k$ differ in the noise-related term. Finally, note that there is an additional recursive equation $\bar{\gamma}_k$ for the additive noise case in contrast with the multiplicative noise case.

\begin{remark}
    The optimal cost  is 
    $f_0(x_0, \bar{x}_0) = {\alpha}_0 var({x}_0) + \bar{\alpha}_0 \bar{x}_0^{2p}$ meaning that, 
as the variance of noises $ \mathbb{E}\left( \epsilon_{k+1}^2  \right)$ vanishes, we retrieve the solution of  the deterministic optimal control. \hfill $\square$
\end{remark}

We have studied how to handle high-order costs that incorporate also the second moment risk term given by the variance. However, the high-order terms enable us to also consider higher moment risk terms. In the following Section, we pursue to minimize higher-order moments.

\section{Higher-Order Moment Costs with  Multiplicative State-Control and Mean-Field Type Dependent Noise}
\label{sec:Multiplicative_mean_field_noise}

We now provide solution for discrete-time optimal control problems with $2o-$moment costs of the form $\mathbb{E}[(x_k-\bar{x}_k)^{2o}]$ and $\mathbb{E}[(u_k-\bar{u}_k)^{2o}]$, where  $o \geq 1. $ The solutions include the optimal control law, the recursive relation for the cost-to-go coefficient $\alpha_k,$ and the expression for the optimal cost. 
%
We consider the following discrete-time system of mean-field type:
\begin{align}
\label{eq:dynamics_4}
   {x}_{k+1} &=  \left(\bar{a}_k \bar{x}_k + \bar{b}_k \bar{u}_k\right) \notag\\
   &\quad\quad\quad\quad +
   \left({a}_k ({x}_k-\bar{x}_k) + {b}_k ({u}_k-\bar{u}_k)\right)\epsilon_{k+1},
\end{align}
with $k \in \{ 0, \dots, N-1\}$, where $x_k \in \mathbb{R}$ is the state, $u_k \in \mathbb{R}$ is the control input, and $\bar{a}_k, \bar{b}_k \in \mathbb{R}$ are system parameters. The objective is to minimize the  $2o$-th moment-dependent  cost function:
\begin{align}
\label{eq:cost_high_moment}
    L(x_0,u) &=  {q}_N {\mathbb{E}[(x_N-\bar{x}_N)^{2o}] }+
    \bar{q}_N \bar{x}_N^{2p} \notag\\
    &\quad\quad\quad\quad + \sum_{k=0}^{N-1} \bigg(  {q}_k {\mathbb{E}[(x_k-\bar{x}_k)^{2o}]} +
    \bar{q}_k \bar{x}_k^{2p} \notag\\
    &\quad\quad\quad\quad + {r}_k {\mathbb{E}[(u_k-\bar{u}_k)^{2o}]}+
     \bar{r}_k \bar{u}_k^{2p}\bigg),
\end{align}
where ${q}_k, {q}_N, {r}_k , \bar{q}_k, \bar{q}_N, \bar{r}_k > 0$  are weights. The control problem is as follows:
\begin{align}
\label{eq:problem4}
    \min_{{u} \in \mathcal{U}_k} \mathbb{E} [L(x_0,{u})],~\mathrm{s.t.}~\eqref{eq:dynamics_4},~x_0~\mathrm{given},
\end{align}
\textcolor{black}{where set of control actions at time step $k$ is 
${U}_k= \{ {u}_k\in \mathbb{R} | \bar{u}_k^{2p} < +\infty \}=\mathbb{R},$ with $\mathcal{U}_k$ as in \eqref{eq:cal_U}.}
\begin{propos}
\label{propos:proposition_7}
    The optimal control input for the higher-order moment problem in \eqref{eq:problem4} with dynamics in \eqref{eq:dynamics_4} and cost functional \eqref{eq:cost_high_moment}, is given by
    \begin{align}
    \label{eq:control_high_moment}
        u_k^* &= - \frac{ \left(\frac{{\alpha}_{k+1}b_km_{k+1, 2o} }{r_k}\right)^{{1}/{(2o-1)}}a_k}{1+ \left(\frac{{\alpha}_{k+1}b_km_{k+1, 2o} }{r_k} \right)^{{1}/{(2o-1)}}b_k}  ({x}_k-\bar{x}_k) \notag\\
        &\quad\quad\quad\quad\quad\quad\quad\quad\, -\frac{\left(\frac{\bar{\alpha}_{k+1} \bar{b}_k}{\bar{r}_k}\right)^{1/(2p-1)} \bar{a}_k }{1 + \left(\frac{\bar{\alpha}_{k+1} \bar{b}_k}{\bar{r}_k}\right)^{1/(2p-1)} \bar{b}_k} \bar{x}_k,
    \end{align}
    where $\bar{\alpha}$ and ${\alpha}$ solve the following recursive equations:
    \begin{subequations}
    \label{eq:recursive_high_moment}
    \begin{align}
    \bar{\alpha}_k &= \bar{q}_k + \bar{r}_k \left(\frac{\left(\frac{\bar{\alpha}_{k+1} \bar{b}_k}{\bar{r}_k}\right)^{1/(2p-1)} \bar{a}_k}{1 + \left(\frac{\bar{\alpha}_{k+1} \bar{b}_k}{\bar{r}_k}\right)^{1/(2p-1)} \bar{b}_k}\right)^{2p} \notag\\
    &+ \bar{\alpha}_{k+1} \left(\bar{a}_{k}- \bar{b}_k \frac{\left(\frac{\bar{\alpha}_{k+1} \bar{b}_k}{\bar{r}_k}\right)^{1/(2p-1)} \bar{a}_k}{1 + \left(\frac{\bar{\alpha}_{k+1} \bar{b}_k}{\bar{r}_k}\right)^{1/(2p-1)} \bar{b}_k}\right)^{2p},
     \\
     %
        {\alpha}_k &= {q}_k + {r}_k \left(\frac{\left(\frac{{\alpha}_{k+1} {b}_k m_{k+1,2o} }{{r}_k}\right)^{1/(2o-1)} a_k}{1 + \left(\frac{{\alpha}_{k+1} {b}_k  m_{k+1,2o}}{{r}_k}\right)^{1/(2o-1)} {b}_k}\right)^{2o}\\ 
        &+ {\alpha}_{k+1} \left(a_k- {b}_k \frac{\left(\frac{{\alpha}_{k+1} {b}_k m_{k+1,2o} }{{r}_k}\right)^{1/(2o-1)} a_k}{1 + \left(\frac{{\alpha}_{k+1} {b}_k m_{k+1,2o}}{{r}_k}\right)^{1/(2o-1)} {b}_k}\right)^{2o} ,\notag
\end{align}
\end{subequations}
with $\bar{\alpha}_N= \bar{q}_N$ and ${\alpha}_N= {q}_N$ as terminal conditions, and $m_{k+1,2o}  = \mathbb{E}[\epsilon_{k+1}^{2o }]$. \hfill $\square$
\end{propos}
Note that the above analysis covers noise beyond the Markovian noise. We now verify the proposed solution by checking the case where the noise process is Markovian.
Comparing the results in Proposition \ref{propos:proposition_6} and Proposition \ref{propos:proposition_7}, we observe the same structure in $\bar{\alpha}_k$. Nevertheless, the multiplicative noise depending on mean-field terms together with the higher-cost functional creates a power dependency for the recursive equation in $\alpha_k$ for Proposition \ref{propos:proposition_7} in contrast to Proposition \ref{propos:proposition_6}.

\subsection{Verification Method in the Space of Measures}

 Let us work directly on the space of probability measures associated with $x_k$ and denote by $\mu_k= \mathbb{P}_{x_k} $ its probability distribution. Then $\mu_{k+1}$ can be obtained  $\mathbb{P}_{\epsilon_{k+1}}$,  $\mu_{k}$ and $x_k,\bar{x}_k,u_k, \bar{u}_k.$ Note here that the mean-field type are involved in the transition kernel.  The cost can also be written in terms of   $\mu_{k}$ and the functional $u_k$ of mean-field type.
Our new augmented state becomes $\mu$  and the system becomes a big deterministic system in the space of measures. The advantage now is that one can write a dynamic programming principle of this new augmented state,
\begin{align*} 
V_k&(\mu_k) = \inf_{u_k} \Bigg(  {q}_k  \int \left[ \left( z-\int y \mu_k(dy) \right)^{2o} \right] \mu_k(dz)\\
&+
    \bar{q}_k \left(\int y \mu_k(dy)\right)^{2p}   \\
    &+ {r}_k \int \left[ \left( u_k(z,\mu_k)-\int  u_k(y, \mu_k) \mu_k(dy) \right)^{2o}\right] \mu_k(dz)\\
    &+ \bar{r}_k  \left(\int u(y,\mu_k) \mu_k(dy)\right)^{2p}  + V_{k+1}(\mu_{k+1})\Bigg)
\end{align*}
and the terminal value is 
\begin{align*}
  V_N(\mu_N)&={q}_N \int \left[\left(z-\int y \mu_N(dy)\right)^{2o}\right] \mu_N(dz) \\
  &+ \bar{q}_N \left(\int y \mu_N(dy)\right)^{2p},  \\
\mu_{k+1}(dy) &= \left\{  
\begin{array}{ll}
  \mathbb{P}_{\epsilon_{k+1} }\bigg( \left({a}_k ({x}_k-\bar{x}_k) + {b}_k ({u}_k-\bar{u}_k)\right)^{-1}\\
  \left(dy-  \left(\bar{a}_k \bar{x}_k + \bar{b}_k \bar{u}_k\right) \right) \bigg) \ \\ 
  \mbox{if} ~ \left({a}_k ({x}_k-\bar{x}_k) + {b}_k ({u}_k-\bar{u}_k) \right)\neq 0,\\\\
\delta_{\bar{a}_k \mathbb{E}{x_k} + \bar{b}_k \mathbb{E}{u_k}} (dy)  \ \mbox{otherwise}.
\end{array} 
\right.
\end{align*}
We seek for $V_k(\mu_k)$ in the form 
\begin{align*}
 \alpha_k  \int \left[\left(z-\int y \mu_k(dy)\right)^{2o}\right] \mu_k(dz) +  \bar{\alpha}_k \left(\int y \mu_k(dy)\right)^{2p}.   
\end{align*}
Inserting the above yields the same dynamics for ${\alpha}_k$ and $\bar{\alpha}_k$ as well as the optimal control structure
\begin{align*}
u_k(x, &\mu_k)= \\
&- \frac{ \left(\frac{{\alpha}_{k+1}b_km_{k+1, 2o} }{r_k}\right)^{{1}/{(2o-1)}}a_k}{1+ \left(\frac{{\alpha}_{k+1}b_km_{k+1, 2o} }{r_k}\right)^{{1}/{(2o-1)}}b_k} \left( {x}-\int y \mu_k(dy)\right) \\
&-\frac{\left(\frac{\bar{\alpha}_{k+1} \bar{b}_k}{\bar{r}_k}\right)^{1/(2p-1)} \bar{a}_k }{1 + \left(\frac{\bar{\alpha}_{k+1} \bar{b}_k}{\bar{r}_k}\right)^{1/(2p-1)} \bar{b}_k}  \left(\int y \mu_k(dy)\right).    
\end{align*} 
Note that we have retrieved the same structure for the optimal control presented in \eqref{eq:control_high_moment} in  Proposition \ref{propos:proposition_7}. We have been addressing optimal control problems, more specifically, we have started from high-order deterministic control problem and have developed it to consider a higher-moment risk-aware stochastic optimal control problems. 
\section{Proofs}
\label{sec:proofs}

\begin{proof}[Proof of Lemma \ref{lema:highorderterms}]
The function $f$ is twice differentiable as it is polynomial. It is strictly convex if its second derivative $f''(z) > 0$, for all $z$.
We differentiate twice to obtain:
\[
f''(z) 
= 2p(2p-1)z^{2p-2} + 2p(2p-1)a^2 (az+b)^{2p-2}
\]
We observe:
\begin{itemize}
\item When $p\geq 1,$ the term $2p(2p-1)z^{2p-2}$ is \textit{positive} for all $z \neq 0$.
\item When $p\geq 1, a, b \neq 0$ The term $2p(2p-1)a^2 (az+b)^{2p-2}$ is always \textit{positive}, for $z \neq - b/a \neq 0$. 
\end{itemize}
Since both terms are non-negative and do not vanish simultaneously, the function is strictly convex under the assumption that $p\geq 1, a, b \neq 0.$ This complete the proof. \qed
\end{proof}

\vspace{0.1cm}
\begin{proof}[Proof of Proposition \ref{propos:problem1}]
    We propose a candidate for the optimal cost  function at time $k$ to be:
\begin{equation}
\label{eq:value_function1}
    f_k(\bar{x}_k) = \bar{{\alpha}}_k \bar{x}_k^4,
\end{equation}
where $\bar{{\alpha}}_k > 0$ is a scalar coefficient to be determined recursively. Then, we use telescopic sum as in \cite{ref1}, i.e.,
\begin{align*}
    f_N(\bar{x}_N) = f_0(\bar{x}_0)+ \sum_{k=0}^{N-1}  f_{k+1}(\bar{x}_{k+1})-f_{k}(\bar{x}_{k})  .
\end{align*}
By taking the difference $L(x_0,u) - f_0(\bar{x}_0)$ we obtain: 
\begin{align*}
   L(x_0,\bar{u})& - f_0(\bar{x}_0) = \bar{q}_N \bar{x}_N^4-f_N(\bar{x}_N) \\
   &+ \sum_{k=0}^{N-1} \left(\bar{q}_k \bar{x}_k^4 + \bar{r}_k \bar{u}_k^4\right)+ f_{k+1}(\bar{x}_{k+1})-f_{k}(\bar{x}_{k}).
\end{align*}
Now, let us introduce the following term:
\begin{align*}
    C_k = \min_{\bar{u}_k} \left( \bar{r}_k \bar{u}_k^4 + \bar{\alpha}_{k+1} (\bar{a}_k \bar{x}_k + \bar{b}_k \bar{u}_k)^4 \right),
\end{align*}
and it yields 
\begin{align*}
   &L(x_0,\bar{u}) - f_0(\bar{x}_0) = \bar{q}_N \bar{x}_N^4-f_N(\bar{x}_N) \\
   &+ \sum_{k=0}^{N-1} \left( \bar{r}_k \bar{u}_k^4+ f_{k+1}(\bar{x}_{k+1})\right)-C_k 
   +  \left(C_k+\bar{q}_k \bar{x}_k^4 -f_{k}(\bar{x}_{k})\right).
\end{align*}

Note here that we are not using the completing the squares technique  but a completion of convex terms. This generalization to convex completion is very important as we will see in the sequel to handle higher-order cost.  The advantage of convex completion is that it covers more general convex functionals that are invariant by   the structure of the control actions.  Here the running cost function is convex in $(\bar{x}_k,\bar{u}_k).$
To find the optimal control, we differentiate the cost function  with respect to $\bar{u}_k:$
\begin{align*}
   4 \bar{r}_k \bar{u}_k^3 + 4 \bar{\alpha}_{k+1} (\bar{a}_k \bar{x}_k + \bar{b}_k \bar{u}_k)^3 \bar{b}_k = 0.
\end{align*}
Simplify and isolate  $(\bar{a}_k \bar{x}_k + \bar{b}_k \bar{u}_k)^3$ yields
\begin{align*}
    \bar{u}_k^3 = -\frac{\bar{\alpha}_{k+1} \bar{b}_k}{\bar{r}_k} (\bar{a}_k \bar{x}_k + \bar{b}_k \bar{u}_k)^3,
\end{align*}
and taking the cube root one obtains the announced optimal control input $\bar{u}_k^*$ as in \eqref{eq:u_optimal1}. Now, let us compute the recursive equation for  $\bar{\alpha}_k$. Substitute $\bar{u}_k^*$ back into the difference  equation and collect terms proportional to $\bar{x}_k^4$ obtaining \eqref{eq:alpha_optimal1}, completing the proof. \qed
\end{proof}

\vspace{0.1cm}
\begin{proof}[Proof of Proposition \ref{propos:problem2}]
The candidate cost is assumed to have the form:
\begin{align*}
    f_k(\bar{x}_k) = \bar{\alpha}_k \bar{x}_k^{2p},
\end{align*}
where $\bar{\alpha}_k > 0$ is a coefficient to be determined recursively, i.e., the optimal cost at each time step is proportional to $\bar{x}_k^{2p}$. The difference term now has
\begin{align*}
 C_k = \min_{\bar{u}_k} \left( \bar{r}_k \bar{u}_k^{2p} + \bar{\alpha}_{k+1} (\bar{a}_k \bar{x}_k + \bar{b}_k \bar{u}_k)^{2p} \right).
\end{align*}
Differentiating with respect to $\bar{u}_k$ gives:
\begin{align*}
   2p \, \bar{r}_k \bar{u}_k^{2p-1} + 2p \, \bar{\alpha}_{k+1} \bar{b}_k (\bar{a}_k \bar{x}_k + \bar{b}_k \bar{u}_k)^{2p-1} = 0.
\end{align*}
Rearranging terms yields
\begin{align*}
    \bar{r}_k \bar{u}_k^{2p-1} + \bar{\alpha}_{k+1} \bar{b}_k (\bar{a}_k \bar{x}_k + \bar{b}_k \bar{u}_k)^{2p-1} = 0.
\end{align*}
Isolate  $(\bar{a}_k \bar{x}_k + \bar{b}_k \bar{u}_k)^{2p-1}$, i.e.,
\begin{align*}
    \bar{u}_k^{2p-1} = -\frac{\bar{\alpha}_{k+1} \bar{b}_k}{\bar{r}_k} (\bar{a}_k \bar{x}_k + \bar{b}_k \bar{u}_k)^{2p-1}.
\end{align*}
Now, by taking the $(2p-1)$-th root on both sides one obtains the optimal control input as in \eqref{eq:u_optimal2}. Finally, by substituting the optimal control $\bar{u}_k^*$ back into the candidate cost yields the recursive equation announced in \eqref{eq:alpha_optimal2}, completing the proof. \qed
\end{proof}

\vspace{0.1cm}
\begin{proof}[Proof of Proposition \ref{propos:proposition_3}]
We propose a candidate cost function at time $k$ to be:
\begin{align*}
    f_k(x_k, \bar{x}_k) =  {\alpha}_k  {var(x_k)}+\bar{\alpha}_k  \bar{x}_k^4+\bar\gamma_k,
\end{align*}
where $ {\alpha}_k, \bar{\alpha}_k, \bar\gamma_k  \geq   0$ are  scalar coefficients to be determined recursively. Now, applying telescopic sum yields    
\begin{align*}
    f_N(x_N) = f_0(x_0)+ \sum_{k=0}^{N-1}  f_{k+1}(x_{k+1})-f_{k}(x_{k})  .
\end{align*}
By taking the difference $L(x_0,u) - f_0(x_0)$ we obtain: 
\begin{align*}
   L(x_0,u)&- f_0(x_0) =  {q}_N {var(x_N)} +\bar{q}_N \bar{x}_N^4-f_N(x_N)\\
   &+ \sum_{k=0}^{N-1} \left( {q}_k {var(x_k) }+\bar{q}_k \bar{x}_k^4 + {r}_k {var(u_k)}+ \bar{r}_k \bar{u}_k^4 \right)\\
   &+ f_{k+1}(x_{k+1})-f_{k}(x_{k}).
\end{align*}
Let us now consider the term
\begin{align*}
    C_k &= \min_{u_k} 
    \big( {r}_k {var(u_k)} +  {\alpha}_{k+1}  {var(x_{k+1})} 
    +\bar{r}_k \bar{u}_k^4 \\
    &+ \bar{\alpha}_{k+1} (\bar{a}_{k}\bar{x}_k + \bar{b}_k \bar{u}_k)^4 \big),
\end{align*}
and
\begin{align}
   L(x_0,u)&- f_0(x_0) =  {q}_N {var(x_N)} +
    \bar{q}_N \bar{x}_N^4-f_N(x_N) \notag\\ 
    &+ \sum_{k=0}^{N-1} \left( {r}_k {var(u_k)}+\bar{r}_k \bar{u}_k^4+ f_{k+1}(x_{k+1})\right)-C_k \notag\\
   &+  \left(C_k+ {q}_k {var(x_k)}+\bar{q}_k \bar{x}_k^4 -f_{k}(x_{k})\right).
   \label{eq:step_3}
\end{align}
Now, the following term can be re-written as follows:
\begin{align*}
{r}_k {var(u_k)} &+  {\alpha}_{k+1}  {var(x_{k+1})} 
= \mathbb{E} ({r}_k {(u_k- \bar{u}_k)^2}  \\
&+  {\alpha}_{k+1}  
( \bar{a}_k (x_k-\bar{x}_k) + \bar{b}_k (u_k-\bar{u}_k)+\epsilon_{k+1} )^2),
\end{align*}
 which means that 
%
%
%
\begin{align*}
 {(u_k- \bar{u}_k)}  = -  \frac{ {\alpha}_{k+1} \bar{b}_k   \bar{a}_k }
{r_k +\alpha_{k+1} \bar{b}^2_k} (x_k-\bar{x}_k).
\end{align*}
To find the deterministic part of optimal control $ \bar{u}_k^*$ , we differentiate the cost function with respect to $u_k$, i.e.,
\begin{align*}
   4 \bar{r}_k \bar{u}_k^3 + 4 \bar{\alpha}_{k+1} (\bar{a}_{k}\bar{x}_k + \bar{b}_k \bar{u}_k)^3 \bar{b}_k = 0.
\end{align*}
Simplify and isolate  $(\bar{a}_{k}\bar{x}_k + \bar{b}_k \bar{u}_k)^3$:
\begin{align*}
    \bar{u}_k^3 = -\frac{\bar{\alpha}_{k+1} \bar{b}_k}{\bar{r}_k} (\bar{a}_{k}\bar{x}_k + \bar{b}_k \bar{u}_k)^3.
\end{align*}
Taking the cube root:
\begin{align*}
    \bar{u}_k^* = -
    \frac{\left(\frac{\bar{\alpha}_{k+1} \bar{b}_k}{\bar{r}_k}\right)^{1/3} \bar{a}_k }{1 + \left(\frac{\bar{\alpha}_{k+1} \bar{b}_k}{\bar{r}_k}\right)^{1/3} \bar{b}_k} \bar{x}_k.
\end{align*}
Therefore, the optimal control is the announced in \eqref{eq:optimal_u_3}. Now, 
substituting $u_k^*$ back into the difference equation in \eqref{eq:step_3} and collecting terms proportional to $\bar{x}_k^4$, one obtains the announced recursive equations for $\bar{\alpha}_k$, $\bar{\gamma}_{k}$, and ${\alpha}_k$, completing the proof. \mbox{~} \qed
\end{proof}

\vspace{0.1cm}
\begin{proof}[Proof of Proposition \ref{propos:proposition_4}]
The candidate cost is assumed to have the form:
\begin{align*}
    f_k(x_k, \bar{x}_k) = {\alpha}_k var({x}_k)+\bar{\alpha}_k \bar{x}_k^{2p} + \bar{\gamma}_k,
\end{align*}
where $\alpha_k, \bar{\alpha}_k  \geq  0$ are coefficients to be determined recursively. The expected value of difference term now has
\begin{align*}
 C_k = \min_{u_k} \bigg(   {r}_k {var(u_k)} +  {\alpha}_{k+1}  &{var(x_{k+1})} +
 \bar{r}_k \bar{u}_k^{2p} \\
 &+ \bar{\alpha}_{k+1} (\bar{a}_{k}\bar{x}_k + \bar{b}_k \bar{u}_k)^{2p} \bigg).
\end{align*}
The stochastic part of the optimal control is determined by 
\begin{align*}
{r}_k {var(u_k)} &+  {\alpha}_{k+1}  {var(x_{k+1})} 
= \mathbb{E} ({r}_k {(u_k- \bar{u}_k)^2}  \\
&+  {\alpha}_{k+1}  
( \bar{a}_k (x_k-\bar{x}_k) + \bar{b}_k (u_k-\bar{u}_k)+\epsilon_{k+1} )^2),
\end{align*}
 which means that 
%
\begin{align*}
    {(u_k- \bar{u}_k)}  = -  \frac{ {\alpha}_{k+1} \bar{b}_k   \bar{a}_k }
{r_k +\alpha_{k+1} \bar{b}^2_k} (x_k-\bar{x}_k).
\end{align*}
The deterministic part of the optimal control is obtained as follows.
Differentiating with respect to $\bar{u}_k$ gives:
\begin{align*}
   2p \, \bar{r}_k \bar{u}_k^{2p-1} + 2p \, \bar{\alpha}_{k+1} \bar{b}_k (\bar{a}_{k}\bar{x}_k + \bar{b}_k \bar{u}_k)^{2p-1} = 0.
\end{align*}
Rearrange terms:
\begin{align*}
    \bar{r}_k \bar{u}_k^{2p-1} + \bar{\alpha}_{k+1} \bar{b}_k (\bar{a}_{k}\bar{x}_k + \bar{b}_k \bar{u}_k)^{2p-1} = 0.
\end{align*}
Isolate  $(\bar{a}_{k}\bar{x}_k + \bar{b}_k \bar{u}_k)^{2p-1}:$
\begin{align*}
    \bar{u}_k^{2p-1} = -\frac{\bar{\alpha}_{k+1} \bar{b}_k}{\bar{r}_k} (\bar{a}_{k}\bar{x}_k + \bar{b}_k \bar{u}_k)^{2p-1}.
\end{align*}
Take the $(2p-1)$-th root on both sides:
\begin{align*}
    \bar{u}_k = -\frac{\left(\frac{\bar{\alpha}_{k+1} \bar{b}_k}{\bar{r}_k}\right)^{1/(2p-1)} \bar{a}_k }{1 + \left(\frac{\bar{\alpha}_{k+1} \bar{b}_k}{\bar{r}_k}\right)^{1/(2p-1)} \bar{b}_k}\bar{x}_k.
\end{align*}
which leads to concluding that the optimal control input is given by the announced expression in \eqref{eq:optimal_u_4}.
Substituting the optimal control $u_k^*$ back into the candidate cost yields the announced recursive equations for $\bar{\alpha}_k$, $\bar{\gamma}_{k}$, and ${\alpha}_k$, where we have used that  
\begin{align*} 
    (X+ \epsilon_{k+1})^{2} = X^{2} +\epsilon_{k+1}^{2} +  2 X \epsilon_{k+1} ,
\end{align*}
completing the proof. \qed
\end{proof}

\vspace{0.1cm}
\begin{proof}[Proof of Proposition \ref{propos:proposition_5}]
We propose a candidate cost  function at time $k$ to be:
\begin{align*}
    f_k(x_k, \bar{x}_k) =  {\alpha}_k  {var(x_k)} +\bar{\alpha}_k  \bar{x}_k^4,
\end{align*}
where ${\alpha}_k, \bar{\alpha}_k  \geq   0$ are  scalar coefficients to be determined recursively. We proceed now with the telescopic sum, i.e.,
\begin{align*}
    f_N(x_N) = f_0(x_0)+ \sum_{k=0}^{N-1}  f_{k+1}(x_{k+1})-f_{k}(x_{k}).
\end{align*}
By taking the difference $L(x_0,u) - f_0(x_0)$ we obtain: 
\begin{align*}
   L(x_0,u) &- f_0(x_0) =  {q}_N {var(x_N)} +\bar{q}_N \bar{x}_N^4-f_N(x_N)\\
   &+ \sum_{k=0}^{N-1} \left( {q}_k {var(x_k)}+\bar{q}_k \bar{x}_k^4 + {r}_k {var(u_k)}+ \bar{r}_k \bar{u}_k^4\right)\\
   &+ f_{k+1}(x_{k+1})-f_{k}(x_{k}). 
\end{align*}
Let us introduce the following term:
\begin{align*}
    C_k = \min_{u_k} 
    \big( {r}_k {var(u_k)} +  {\alpha}_{k+1} & {var(x_{k+1})} 
    +\bar{r}_k \bar{u}_k^4 \\
    &+ \bar{\alpha}_{k+1} (\bar{a}_{k}\bar{x}_k + \bar{b}_k \bar{u}_k)^4 \big).
\end{align*}
Then, the difference $L(x_0,u) - f_0(x_0)$ is re-written as
\begin{align*}
   L(x_0,u)&- f_0(x_0) =  {q}_N {var(x_N)} +
    \bar{q}_N \bar{x}_N^4-f_N(x_N) \\ 
    &+ \sum_{k=0}^{N-1} \left( {r}_k {var(u_k)}+\bar{r}_k \bar{u}_k^4+ f_{k+1}(x_{k+1})\right)\\
    &-C_k 
   +  \left(C_k+ {q}_k {var(x_k)} +\bar{q}_k \bar{x}_k^4 -f_{k}(x_{k})\right),
\end{align*}
where
\begin{align*}
{r}_k &{var(u_k)} +  {\alpha}_{k+1}  {var(x_{k+1})} 
= \mathbb{E} ({r}_k {(u_k- \bar{u}_k)^2}  \\
&+  {\alpha}_{k+1}  
( \bar{a}_k (x_k-\bar{x}_k) + \bar{b}_k (u_k-\bar{u}_k)+(x_k-\bar{x}_k)\epsilon_{k+1}  )^2),
\end{align*}
 which means that 
 %
%
\begin{align*}
{(u_k- \bar{u}_k)}  = -  \frac{ {\alpha}_{k+1} \bar{b}_k   \bar{a}_k }
{r_k +\alpha_{k+1} \bar{b}^2_k} (x_k-\bar{x}_k).   
\end{align*}
To find the deterministic part of optimal control $ \bar{u}_k^*$ , we differentiate the cost function with respect to $\bar{u}_k$, i.e.,
\begin{align*}
   4 \bar{r}_k \bar{u}_k^3 + 4 \bar{\alpha}_{k+1} (\bar{a}_{k}\bar{x}_k + \bar{b}_k \bar{u}_k)^3 \bar{b}_k = 0.
\end{align*}
Simplify and isolate  $(\bar{a}_{k}\bar{x}_k + \bar{b}_k \bar{u}_k)^3$ resulting in
\begin{align*}
    \bar{u}_k^3 = -\frac{\bar{\alpha}_{k+1} \bar{b}_k}{\bar{r}_k} (\bar{a}_{k}\bar{x}_k + \bar{b}_k \bar{u}_k)^3.
\end{align*}
Taking the cube root:
\begin{align*}
    \bar{u}_k^* = -
    \frac{\left(\frac{\bar{\alpha}_{k+1} \bar{b}_k}{\bar{r}_k}\right)^{1/3} \bar{a}_k }{1 + \left(\frac{\bar{\alpha}_{k+1} \bar{b}_k}{\bar{r}_k}\right)^{1/3} \bar{b}_k} \bar{x}_k.
\end{align*}
Then, the optimal control is obtained corresponding to \eqref{eq:optimal_u_5}.
Now, by substituting $u_k^*$ back into the difference equation and collecting terms proportional to $\bar{x}_k^4$, one obtains the recursive equation for  $\alpha_k$, and $\bar{\alpha}_k$ as in \eqref{eq:resursive_eqs}, completing the proof. \qed
\end{proof}

\vspace{0.1cm}
\begin{proof}[Proof of Proposition \ref{propos:proposition_6}]
The candidate cost is assumed to have the form:
\begin{align*}
    f_k(x_k, \bar{x}_k) = {\alpha}_k var({x}_k)+\bar{\alpha}_k \bar{x}_k^{2p} ,
\end{align*}
where $\alpha_k, \bar{\alpha}_k  \geq  0$ are coefficients to be determined recursively.    
The expected value of difference term now has
\begin{align*}
 C_k = \min_{u_k} \bigg(   {r}_k {var(u_k)} +  {\alpha}_{k+1} & {var(x_{k+1})} +
 \bar{r}_k \bar{u}_k^{2p} \\
 &+ \bar{\alpha}_{k+1} (\bar{a}_{k}\bar{x}_k + \bar{b}_k \bar{u}_k)^{2p} \bigg).
\end{align*}
The stochastic part of the optimal control is determined by 
\begin{align*}
{r}_k &{var(u_k)} +  {\alpha}_{k+1}  {var(x_{k+1})} 
= \mathbb{E} ({r}_k {(u_k- \bar{u}_k)^2}  \\
&+  {\alpha}_{k+1}  
( \bar{a}_k (x_k-\bar{x}_k) + \bar{b}_k (u_k-\bar{u}_k)+(x_k-\bar{x}_k)\epsilon_{k+1} )^2),
\end{align*}
 which means that 
\begin{align*}
{(u_k- \bar{u}_k)}  = -  \frac{ {\alpha}_{k+1} \bar{b}_k   \bar{a}_k }
{r_k +\alpha_{k+1} \bar{b}^2_k} (x_k-\bar{x}_k).    
\end{align*}
The deterministic part of the optimal control is obtained as follows.
Differentiating with respect to $\bar{u}_k$ gives:
\begin{align*}
   2p \, \bar{r}_k \bar{u}_k^{2p-1} + 2p \, \bar{\alpha}_{k+1} \bar{b}_k (\bar{a}_{k}\bar{x}_k + \bar{b}_k \bar{u}_k)^{2p-1} = 0.
\end{align*}
Rearrange terms:
\begin{align*}
    \bar{r}_k \bar{u}_k^{2p-1} + \bar{\alpha}_{k+1} \bar{b}_k (\bar{a}_{k}\bar{x}_k + \bar{b}_k \bar{u}_k)^{2p-1} = 0.
\end{align*}
Isolate  $(\bar{a}_{k}\bar{x}_k + \bar{b}_k \bar{u}_k)^{2p-1}:$
\begin{align*}
    \bar{u}_k^{2p-1} = -\frac{\bar{\alpha}_{k+1} \bar{b}_k}{\bar{r}_k} (\bar{a}_{k}\bar{x}_k + \bar{b}_k \bar{u}_k)^{2p-1}.
\end{align*}
Take the $(2p-1)$-th root on both sides:
\begin{align*}
    \bar{u}_k = -\frac{\left(\frac{\bar{\alpha}_{k+1} \bar{b}_k}{\bar{r}_k}\right)^{1/(2p-1)} \bar{a}_k }{1 + \left(\frac{\bar{\alpha}_{k+1} \bar{b}_k}{\bar{r}_k}\right)^{1/(2p-1)} \bar{b}_k} \bar{x}_k,
\end{align*}
and we obtain the announced optimal control input ${u}^*_k$ in \eqref{eq:optimal_u_6}. 
Substituting the optimal control $u_k^*$ back into the candidate cost yields the recursive equations for ${\alpha}_k$, and $\bar{\alpha}_{k}$  in \eqref{eq:resursive_eqs_2}, where we have used that  
\begin{align*}
    (a Y+ Y\epsilon_{k+1})^{2} = Y^2(a^{2} +\epsilon_{k+1}^{2} +  2 a\epsilon_{k+1}),
\end{align*}
completing the proof. \qed
\end{proof}

\vspace{0.1cm}
\begin{proof}[Proof of Proposition \ref{propos:proposition_7}]
    The candidate cost is assumed to have the form:
\begin{align*}
    f_k(x_k, \bar{x}_k) =  {\alpha}_{k} {\mathbb{E}[(x_k-\bar{x}_k)^{2o}]} +\bar{\alpha}_k \bar{x}_k^{2p},
\end{align*}
where $\alpha_k \geq  0$,  $\bar\alpha_k\geq  0$ are coefficients to be determined recursively.
The expected value of difference term now has
\begin{align*}
 C &= \min_{u_k} \bigg(   {r}_k {\mathbb{E}[(u_k-\bar{u}_k)^{2o}] } +  {\alpha}_{k+1}  {\mathbb{E}[(x_{k+1}-\bar{x}_{k+1})^{2o}]}\\
 &+\bar{r}_k \bar{u}_k^{2p} + \bar{\alpha}_{k+1} (\bar{a}_{k}\bar{x}_k + \bar{b}_k \bar{u}_k \bigg)^{2p} 
   \bigg).
\end{align*}
This minimization is subject to $\mathcal{F}_k$ adaptivity of optimal control. 
We compute 
\begin{align*}
\mathbb{E}&[(x_{k+1}-\bar{x}_{k+1})^{2o}] \\
&=   
\mathbb{E}[({a}_k ({x}_k-\bar{x}_k) + {b}_k ({u}_k-\bar{u}_k))\epsilon_{k+1})^{2o}],\\ 
&=  \mathbb{E}[({a}_k ({x}_k-\bar{x}_k) + {b}_k ({u}_k-\bar{u}_k))^{2o}]  m_{k+1, 2o}.
\end{align*}
\begin{figure}[t!]
    \centering
    \begin{tabular}{cc}
       \includegraphics[width=0.5\columnwidth]{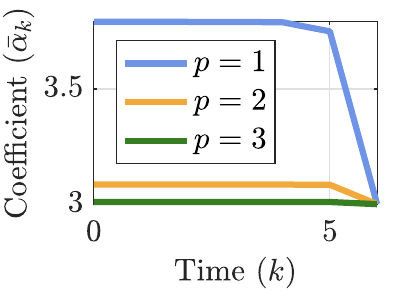}  &
       \hspace{-0.5cm}
       \includegraphics[width=0.5\columnwidth]{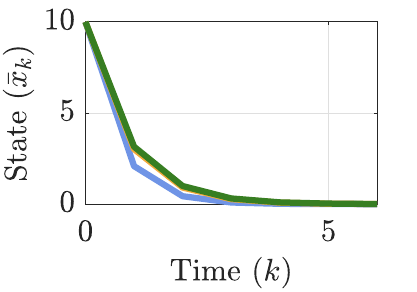}  \\
    \end{tabular}
    \caption{Results example corresponding to Proposition \ref{propos:problem1} and Proposition \ref{propos:problem2}.}
    \label{fig:p1_2}
\end{figure}
The stochastic part of the control yields
\begin{align*}
 u_k-\bar{u}_k= - \frac{ \left(\frac{{\alpha}_{k+1}b_km_{k+1, 2o} }{r_k}\right)^{{1}/{(2o-1)}}a_k}{1+ \left(\frac{{\alpha}_{k+1}b_km_{k+1, 2o} }{r_k}\right)^{{1}/{(2o-1)}}b_k}  ({x}_k-\bar{x}_k).
\end{align*}
The deterministic part of the control yields
\begin{align*}
    \bar{u}_k = -\frac{\left(\frac{\bar{\alpha}_{k+1} \bar{b}_k}{\bar{r}_k}\right)^{1/(2p-1)} \bar{a}_k}{1 + \left(\frac{\bar{\alpha}_{k+1} \bar{b}_k}{\bar{r}_k}\right)^{1/(2p-1)} \bar{b}_k} \bar{x}_k.
\end{align*}
We can easily check that the obtained optimal control is the announced in \eqref{eq:control_high_moment}, which  is adapted to the filtration generated by $\{ x_0, \ldots, x_k\}$. 
Substituting the optimal control $u_k^*$ back into the candidate cost yields \eqref{eq:recursive_high_moment}, and the optimal cost is 
\begin{align*}
    f_0(x_0, \bar{x}_0) = {\alpha}_0 {\mathbb{E}[(x_0-\bar{x}_0)^{2o}]}+ \bar{\alpha}_0 \bar{x}_0^{2p},
\end{align*}
completing the proof. \qed
\end{proof}

\section{Numerical Examples}
\label{sec:examples}

\begin{figure}[t!]
    \centering
    \begin{tabular}{cc}
       \includegraphics[width=0.5\columnwidth]{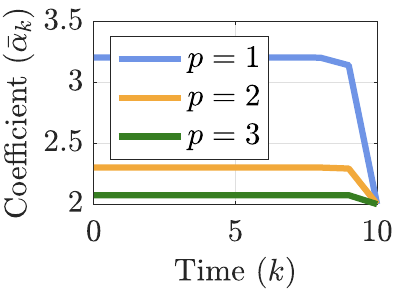}  &
       \hspace{-0.5cm}
       \includegraphics[width=0.5\columnwidth]{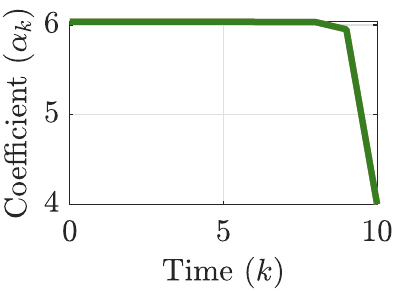}  \\
         \includegraphics[width=0.5\columnwidth]{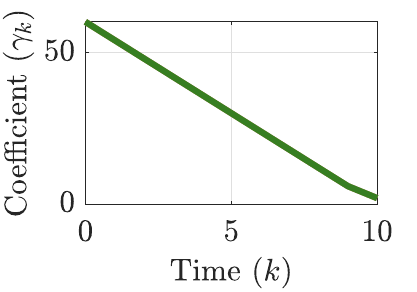}  &
       \hspace{-0.5cm}
       \includegraphics[width=0.5\columnwidth]{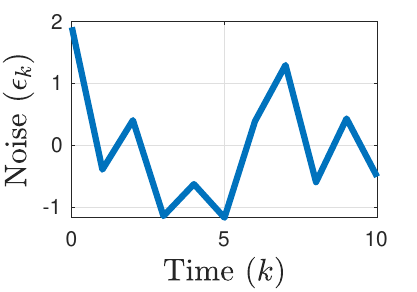}  \\
       \includegraphics[width=0.5\columnwidth]{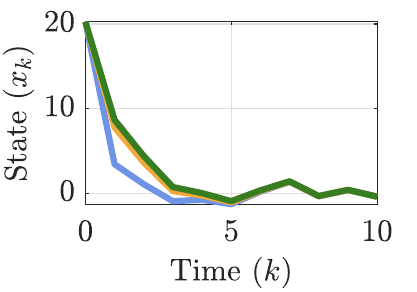}  &
       \hspace{-0.5cm}
       \includegraphics[width=0.5\columnwidth]{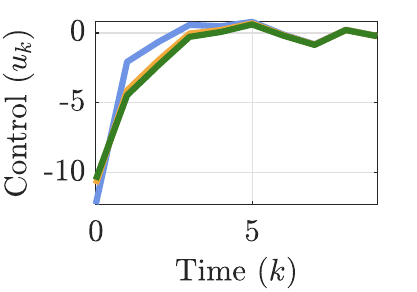}  \\
    \end{tabular}
    \caption{Results example corresponding to Proposition \ref{propos:proposition_3} and Proposition \ref{propos:proposition_4}.}
    \label{fig:p3_4}
    
       \vspace{-0.5cm}
\end{figure}

In this Section we show some numerical examples that allow illustrating the theoretical results presented in this paper and corresponding to each one of the Propositions.

\vspace{-0.25cm}
\subsection*{Example 1: Deterministic control with higher-order costs} 

This example corresponds to the results presented in Proposition \ref{propos:problem1} and Proposition \ref{propos:problem2}. We study the effect of $p$ on $\bar{\alpha}$, $x$, and $\bar{u}$. Let us consider a system with $\bar{a}_k = 1$, $\bar{b}_k = 2$, $\bar{q}_k = \bar{q}_N = 3$, $\bar{r}_k = 4$ and $\bar{x}_0 = 10$; and the terminal time $N=6$. The results are presented in Figure \ref{fig:p1_2}. It can be seen the backward recursive equation evolution $\bar{\alpha}$ depending on the parameter $p$. For all the cases, and as expected, we observe that the state goes to zero minimizing the cost. Note that, according to the optimal control action in Proposition \ref{propos:problem2}, and given that $\bar{\alpha}_k$ is decreasing, and $\bar{x}_k$ goes to zero, $u^*_k$ starts at a negative value and vanishes along the time as the state goes to zero.

\vspace{-0.25cm}
\subsection*{Example 2: Variance-aware control with additive noise} 

This example corresponds to the results presented in Proposition \ref{propos:proposition_3} and Proposition \ref{propos:proposition_4}. We study the effect of $p$ on $\bar{\alpha}$, ${\alpha}$, $\gamma$, $x$, and ${u}$. Let us consider a system with $\bar{a}_k = 2$, $\bar{b}_k = 3$, $\epsilon$ as in Figure \ref{fig:p3_4} with unitary variance $\sigma=1$, and $\bar{q}_k = \bar{q}_N = 2$, $\bar{r}_k = 3$, ${q}_k = {q}_N = 4$, ${r}_k = 5$ and $\bar{x}_0 = x_0 = 20.25$; and the terminal time $N=10$. The results are presented in Figure \ref{fig:p3_4}. We can observe that the recursive equations for $\alpha_k$ and $\bar{\gamma}_k$ are the same for all the different $p$ values as shown in Proposition \ref{propos:proposition_4}. Figure \ref{fig:p3_4} also shows the evolution of $\bar{\alpha}_k$ exhibiting a decreasing behavior. Regarding the system state, we observe convergence around zero showing that the expected state vanishes along the time minimizing the cost. 

\vspace{-0.25cm}
\subsection*{Example 3: Variance-Aware Stochastic Cost with Multiplicative Noise} 

\begin{figure}[t!]
    \centering
    \begin{tabular}{cc}
       \includegraphics[width=0.5\columnwidth]{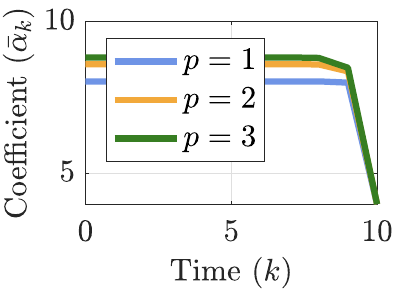}  &
       \hspace{-0.5cm}
       \includegraphics[width=0.5\columnwidth]{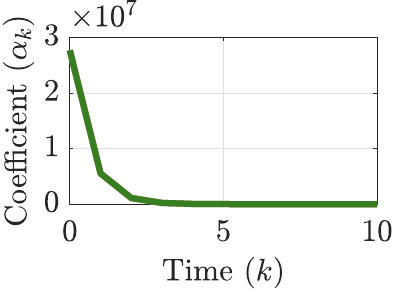}  \\
         \includegraphics[width=0.5\columnwidth]{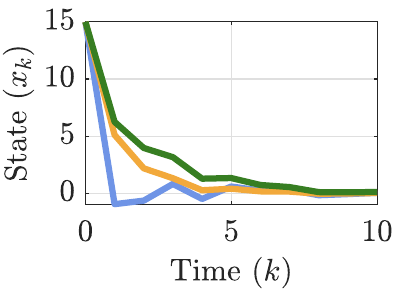}  &
       \hspace{-0.5cm}
       \includegraphics[width=0.5\columnwidth]{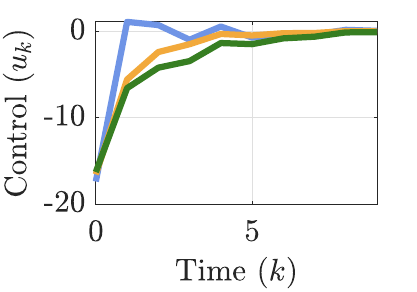}  \\
    \end{tabular}
    \caption{Results example corresponding to Proposition \ref{propos:proposition_5} and Proposition \ref{propos:proposition_6}. The noise $\epsilon$ for this example is the same one presented in Figure \ref{fig:p3_4}.}
    \label{fig:p5_6}
    
       \vspace{-0.5cm}
\end{figure}

This example corresponds to the results presented in Proposition \ref{propos:proposition_5} and Proposition \ref{propos:proposition_6}. We study the effect of $p$ on $\bar{\alpha}$, ${\alpha}$, $x$, and ${u}$. Let us consider a system with $\bar{a}_k = 7$, $\bar{b}_k = 6$, $\epsilon$ as in Figure \ref{fig:p3_4} with unitary variance $\sigma=1$, and $\bar{q}_k = \bar{q}_N = 4$, $\bar{r}_k = 3$, ${q}_k = {q}_N = 2$, ${r}_k = 1$ and $\bar{x}_0 =15.3$, $x_0 =15$; and the terminal time $N=10$. The results are presented in Figure \ref{fig:p5_6}. We can observe that the recursive solution for $\alpha_k$ across the different $p$ values are the same according to the result in Proposition \ref{propos:proposition_6}. Figure \ref{fig:p5_6} also shows the decreasing behavior of $\bar{\alpha}_k$ for the different $p$ values. Also, the state evolves converging around zero, i.e., the expected value of the state vanishes as desired, minimizing the cost. 

\begin{figure}[t!]
    \centering
    \begin{tabular}{cc}
       \includegraphics[width=0.5\columnwidth]{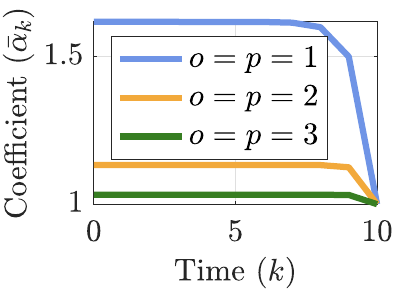}  &
       \hspace{-0.5cm}
       \includegraphics[width=0.5\columnwidth]{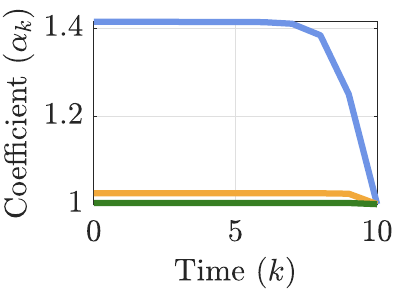}  \\
         \includegraphics[width=0.5\columnwidth]{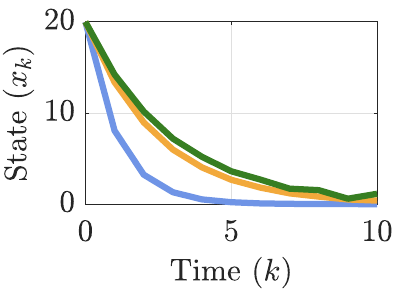}  &
       \hspace{-0.5cm}
       \includegraphics[width=0.5\columnwidth]{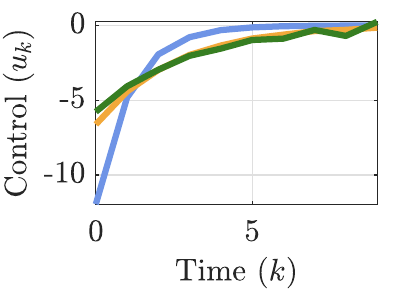}  \\
    \end{tabular}
    \caption{Results example corresponding to Proposition \ref{propos:proposition_7}. The noise $\epsilon$ for this example is the same one presented in Figure \ref{fig:p3_4}.}
    \label{fig:p7}
    
    \vspace{-0.5cm}
\end{figure}

\vspace{-0.25cm}
\subsection*{Example 4: Higher-Order Moment Costs with  Multiplicative State-Control and Mean-Field Type Dependent Noise} 

This example corresponds to the results presented in Proposition \ref{propos:proposition_7}. We study the effect of $p$ on $\bar{\alpha}$, ${\alpha}$, $x$, and ${u}$. Let us consider a system with $\bar{a}_k = \bar{b}_k = 1$, ${a}_k = {b}_k = 1/2$,  $\epsilon$ as in Figure \ref{fig:p3_4} with unitary variance $\sigma=1$, and $\bar{q}_k = \bar{q}_N = \bar{r}_k = 1$, ${q}_k = {q}_N = {r}_k = 1$ and $\bar{x}_0 =20.01$, $x_0 =20$; and the terminal time $N=10$. The results are presented in Figure \ref{fig:p7}. The comparison of the higher-order parameters for $o=p$ are shown and for all the cases the expected value for the system state $x$ decreases to zero as expected. In addition, it is worth noting that the convergence is slower as the power value is increased in the cost functional.

\vspace{-0.25cm}
\subsection*{The Risk Awareness Value}
\begin{figure}[t!]
    \centering
    \begin{tabular}{cc}
       \includegraphics[width=0.5\columnwidth]{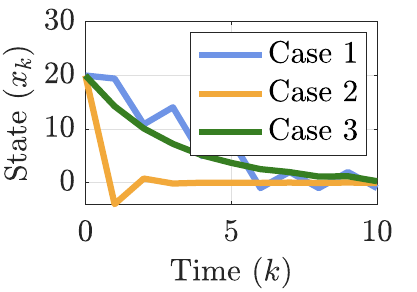}  &
       \hspace{-0.5cm}
       \includegraphics[width=0.5\columnwidth]{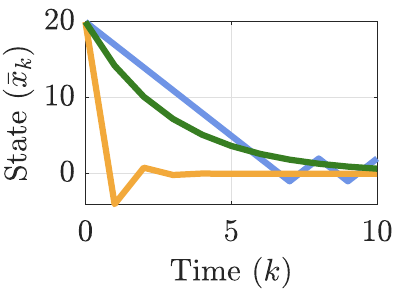}  \\
         \includegraphics[width=0.5\columnwidth]{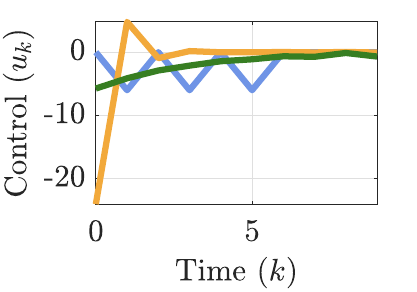}  &
       \hspace{-0.5cm}
       \includegraphics[width=0.5\columnwidth]{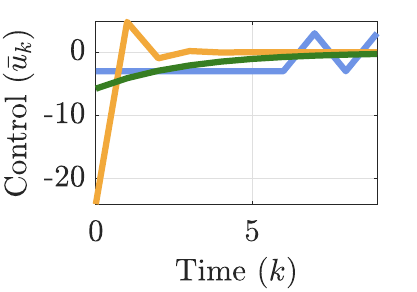}  \\
    \end{tabular}
    \caption{Comparison results: Case 1: simple controller with finite-set control actions, Case 2: state-feedback controller, and Case 3: risk-aware controller.}
    \label{fig:comparison}
    
    \vspace{-0.5cm}
\end{figure}

\begin{table*}[t!]
    \centering
    \caption{Comparison among, Case 1: simple controller with finite-set control actions, Case 2: state-feedback controller, and Case 3: risk-aware controller (see Figure \ref{fig:comparison}).}
    \label{tab:comparison}
    \begin{tabular}{c|ccc|ccc|ccc}
    \hline
    \textbf{Control} & \multicolumn{3}{c|}{$\zeta=1$}& \multicolumn{3}{c|}{$\zeta=2$}& \multicolumn{3}{c}{$\zeta=3$}\\
       \textbf{Strategy}  &  $\mathrm{KPI}^1_x$ & $\mathrm{KPI}^1_u$ & $\mathrm{KPI}^1_x + \mathrm{KPI}^1_u$ &  $\mathrm{KPI}^2_x$ & $\mathrm{KPI}^2_u$ & $\mathrm{KPI}^2_x + \mathrm{KPI}^2_u$ &  $\mathrm{KPI}^3_x$ & $\mathrm{KPI}^3_u$ & $\mathrm{KPI}^3_x + \mathrm{KPI}^3_u$\\
       \hline
       \hline
       Case 1  & 1.19E+03	&180	&1.37E+03	&3.03E+05	&1620	&3.04E+05	&9.80E+07	&14580	&9.81E+07\\
       Case 2 & 417.0897	&600.6092	&1.02E+03	&1.61E+05	&3.33E+05	&4.94E+05	&6.42E+07	&1.92E+08	&2.56E+08\\
       Case 3 & 811.8863	&67.5081	&\textbf{8.79E+02}	&2.16E+05	&1.48E+03	&\textbf{2.17E+05}	&7.38E+07	&4.19E+04	&\textbf{7.39E+07}\\
       \hline
    \end{tabular}
\end{table*}
This Section is devoted to analyze the usefulness of the incorporation of risk awareness to the control design. To appreciate this advantage, we compare a simple controller in which only a finite set of control actions are possible, a state-feedback controller, and a risk-aware higher-order cost controller as the ones we have been studying throughout this paper. We define a unified risk framework to measure the performance of the controllers, i.e., we introduce a risk-aware Key Performance Indicator (KPI) for both terms in $x$ and $u$:
\begin{align*}
    \mathrm{KPI}^\zeta_x &=  \sum_{k=0}^{N} (x_k-\bar{x}_k)^{2\zeta} + \bar{x}_k^{2\zeta},\\
    \mathrm{KPI}^\zeta_u &=  \sum_{k=0}^{N} (u_k-\bar{u}_k)^{2\zeta} + \bar{u}_k^{2\zeta},
\end{align*}
for $\zeta \in \{1,2,3\}$. 
Let us consider the same system dynamics as in Example 4, with the same noise as in Figure \ref{fig:p3_4}, and controlled under three strategies as follows:
\begin{itemize}
    \item \textit{Case 1:} We consider a simple controller with only five feasible control actions, i.e., 
    \begin{align*}
        u^*_k = -3 \mathrm{sgn}(\bar{x}_k- x_k) -3 \mathrm{sgn}(\bar{x}_k),
    \end{align*}
    meaning that $u^*_k \in \{-6,-3,0,3,6\}$, for all $k$.
    \item \textit{Case 2:} Let us consider a more sophisticated controller in comparison to the one in Case 1. Here, we have a state-feedback controller as follows:
    \begin{align*}
        u^*_k = -3(\bar{x}_k- x_k) -3(\bar{x}_k),
    \end{align*}
    meaning that $u^*_k \in \mathbb{R}$ having more flexibility.
    \item \textit{Case 3:} We consider the controller presented in Example 4, i.e., a risk-aware controller with higher-order cost functional with value $o=p=3$.
\end{itemize}
Figure \ref{fig:comparison} presents the comparison of the performance of the three controllers described above and corresponding to Cases 1, 2 and 3. We can observe that all of the controllers accomplish to drive the system state to zero; however, the KPI shows a much better performance for the risk-aware controller (see Table \ref{tab:comparison}). Therefore, we can appreciate the advantage of the risk-awareness in the design of optimal controllers considering higher-order terms.
\section{Discussions and Concluding Remarks}
\label{sec:conclusions}

We  investigated higher-order costs in the context of discrete time mean-field-type optimal control theory, addressing limitations in traditional quadratic cost models and offering solutions for various scenarios.  We presented semi-explicit solutions for different problem types, including deterministic control with higher-order costs, variance-aware control with additive noise, and risk-aware stochastic control with multiplicative noise.  These solutions provide a valuable tool for understanding and controlling a variety of systems that involve non-linear penalties and risk-sensitive behavior.


%

\bibliographystyle{unsrt}
\bibliography{references.bib}
 \balance

\end{document}